\documentclass[11pt]{article}
\usepackage{amssymb}
\usepackage{amsfonts}

\newtheorem{myDef}{Definition}[section]

\newtheorem{myTheorem}[myDef]{Theorem}
\newtheorem{myLemma}[myDef]{Lemma}

\newtheorem{myCorollary}[myDef]{Corollary}

\newenvironment{example}{\noindent{\bf Example}\hspace*{1em}}{\bigskip}
\newenvironment{proof}{\noindent{\bf Proof}\hspace*{1em}}{\qed \bigskip}
\newenvironment{proofof}[1]{\noindent{\bf Proof of #1:}
\hspace*{1em}}{\qed\bigskip}
\newcommand{\qed}{\rule{7pt}{7pt}}

\newcommand{\calJ}{{\cal J}}

\begin{document}

\title{\textbf{Poincar\'{e} type inequalities on the discrete cube \ and in the CAR
algebra}}
\author{L. Ben Efraim\thanks{
The Hebrew University, Jerusalem, Israel, and Weizmann Institute of Science,
Rehovot, Israel. Limor\_be@cs.huji.ac.il. } \ and \ \ F. Lust-Piquard\thanks{
Universit\'{e} de Cergy, France. francoise.piquard@math.u-cergy.fr.} }
\maketitle

\begin{abstract}
We prove $L^{p}$ Poincar\'{e} inequalities with suitable dimension free
constants for functions on $\ $the cube $\{-1,1\}^{n}.$ As well known, such
inequalities for $p$ an even integer allow to recover an exponential
inequality hence the concentration phenomenon first obtained by Bobkov and G%
\"{o}tze. We also get inequalities between the $L^{p}$ norms of $\left\vert
\nabla f\right\vert $ and $\Delta ^{\alpha }f,\alpha >0;$ moreover $L^{p}$
spaces may be replaced by more general ones.

Similar results hold true, replacing functions on the cube by matrices in
the * algebra spanned by $n$ fermions and the $L^{p}$ norm by the Schatten
norm $C_{p}.$
\end{abstract}

\setlength{\baselineskip}{16pt}

\section{\textbf{\ Introduction and some notation}}

Let $\Omega _{n}=\{-1,1\}^{n}$ be the n-dimensional cube (or dyadic group),
equipped with its uniform probability $\mathbb{P}_{n}$ and the corresponding
expectation $\mathbb{E}$. We prove (theorem \ref{theorem:main}) $L^{p}$ $%
(1\leq p<\infty )$ Poincar\'{e} inequalities with suitable constants for the
length of the discrete gradient $\left\vert \nabla f\right\vert ,$where $f$
is a function on $\Omega _{n}.$ As well known, such inequalities for $p$ an
even integer imply an exponential inequality (corollary \ref{cor:exp}) hence
the concentration phenomenon (corollary 4.2), with dimension free constants$%
. $ The exponential inequality was first obtained by Bobkov and G\"{o}tze
\cite{BG}. Their proof relies on a modified Log-Sobolev inequality for
functions on $\Omega _{n},$ which is not involved in ours.

Actually we also get a strengthening of Poincar\'{e} inequalities in two
directions: $f-\mathbb{E}f$ may be replaced by $\Delta ^{\alpha }f,$ $%
0<\alpha <\frac{1}{2},$ where $\Delta $ is the discrete Laplacian, and $%
L^{p} $ norms may be replaced by more general ones, see theorem \ref%
{theorem:main}. \ In theorem \ref{5.1} we prove reverse inequalities
involving $\left\vert \nabla f\right\vert $ and $\Delta ^{\alpha }f,$ $%
\alpha >\frac{1}{2}.$

Using similar ideas, we also get (theorem \ref{theorem:mainCAR}
Poincar\'{e} inequalities for the CAR *-algebra spanned by $n$
fermions, equipped with the Schatten norm $C_{p},$ hence an analogue
of the concentration phenomenon (corollary \ref{cor:concCAR}).
Generalizations and reverse inequalities also hold, as in the cube
case (theorems \ref{theorem:mainCAR} and \ref{6.8})

The main ideas used in the present paper already appeared in \cite{LP} where
Riesz transforms on $L^{p}(\Omega _{n})$ ($1<p<\infty ),$ i.e. the behavior
of $\Vert |\nabla \Delta ^{-\frac{1}{2}}(f)|\Vert _{L^{p}(\Omega _{n})},$
were considered.

The proofs of Theorems \ref{theorem:main} and \ref{theorem:mainCAR} are
similar in spirit to Maurey and Pisier's in the gaussian setting \cite[
theorem 2.2]{P}. As in \cite{LP} the proofs of theorems \ref{5.1} and \ref%
{6.8} are inspired by \cite{P''}. Namely the iid gaussian variables $%
X_{1},..,X_{n}$ are replaced either by iid Bernoulli variables $\varepsilon
_{1},..,\varepsilon _{n}$ (in other words the coordinate functions $\omega
_{1},..,\omega _{n}$ defined on $\Omega _{n}),$ or by fermions $%
Q_{1}^{\prime },..,Q_{n}^{\prime }$. Let us explain how we replace the
independant copy $(Y_{1},..,Y_{n})$ of the gaussian vector $%
(X_{1},..,X_{n}), $ for example in the cube case:

We consider the GNS representation of $L^{\infty }(\Omega _{n})$ in $%
B(L^{2}(\Omega _{n})).$ $L^{2}(\Omega _{n})$ being equipped with the o.n
basis of Walsh functions $\omega _{A},A\subset \{1,.,n\}$, $B(L^{2}(\Omega
_{n})$ is identified with the algebra $\mathcal{M}_{2^{n}}$ of $2^{n}\times
2^{n}$ matrices, which is equipped with its normalized trace; the coordinate
functions $\omega _{1},..,\omega _{n}$ (acting by pointwise multiplication
on $L^{2}(\Omega _{n}))$ are represented by matrices $Q_{1},..,Q_{n}\in
\mathcal{M}_{2^{n}}$. We introduce a sequence of matrices $(P_{1},..,P_{n})$
with properties similar to those of $(Q_{1},..,Q_{n})$, which play the role
of "an independent copy" of $(Q_{1},..,Q_{n})$. The exact definition of $%
(P_{1},..,P_{n})$ and $(Q_{1},..,Q_{n})$ appears in the preliminaries. The
matrices $P_{j}\in \mathcal{M}_{2^{n}}$ have particular commutation
relations with the $Q_{j}$'s.

Let us recall the definitions of $\nabla $ and $\Delta .$ Let $%
e_{j}=(1,..,-1,..,1)\in \Omega _{n},$ where $-1$ occurs at coordinate $j$
and let $\nabla f=(\partial _{j}f)_{j=1}^{n},$ where, for $f\in L^{\infty
}(\Omega _{n}),$

\[
(\partial _{j}f)(x)=f(x)-f(xe_{j}).
\]

\noindent Here$\ xy$ denotes the coordinates product of $x,y\in \Omega _{n}$.

\noindent We denote by $\omega _{j},1\leq j\leq n,$ the $j^{th}$ coordinate
function on $\Omega _{n}$ and, for a non empty subset $A\subset \{1,..,n\},$

\[
\omega _{A}=\prod\limits_{j\in A}\omega _{j},
\]

\noindent while $\mathbf{\omega _{\varnothing}=1}$ is the unit of $L^{\infty
}(\Omega _{n})$. Then

\[
\partial _{j}\omega _{A}=\left\{
\begin{array}{ll}
0 & \mbox{if $j \not\in A$} \\
2\omega _{A} & \mbox{if $j \in A$}%
\end{array}%
\right.
\]

\noindent Actually it will be more convenient to replace $\partial _{j}$ by
another operator

\[
D_{j}=\frac{1}{2}\omega _{j}\partial _{j}
\]

\noindent so that

\[
\left\vert \nabla f\right\vert =(\sum_{j=1}^{n}\left\vert (\partial
_{j}f)(x)\right\vert ^{2})^{\frac{1}{2}}=2(\sum_{j=1}^{n}\left\vert
D_{j}(f)\right\vert ^{2})^{\frac{1}{2}}.
\]

\noindent The discrete Laplacian is $\Delta =\sum_{j=1}^{n}\partial
_{j}^{\ast }\partial _{j}=\sum_{j=1}^{n}\partial _{j}^{2}$ (where $\partial
_{j}^{\ast }$ is the adjoint of $\partial _{j}$ acting on $L^{2}(\Omega
_{n}) $) so that

\[
\Delta (\omega _{A})=4\left\vert A\right\vert \omega _{A}.
\]

\noindent Actually we shall rather consider the number operator

\[
N=\sum_{j=1}^{n}D_{j}^{\ast }D_{j}=\frac{\Delta }{4}.
\]

\noindent A function $f$ on $\Omega _{n}$ can be written as
\[
f=\mathbb{E}f+\sum_{A\subset \{1,..,n\},A\neq \varnothing }\hat{f}%
(A)\,\omega _{A}.\,
\]%
For $\,0\leq \theta <\frac{\pi }{2}$, we have%
\begin{equation}
\cos ^{\frac{\Delta }{4}}\theta (f)\,=\cos ^{N}\theta (f)\,=\mathbb{E}%
f+\,\sum_{A\neq \varnothing }\cos ^{|A|}\theta \,\hat{f}(A)\,\omega _{A}.
\label{+}
\end{equation}

We denote by $E$ a symmetric function space on $\Omega _{n}=\{-1,1\}^{n}$,
i.e. $E$ is a Banach lattice of complex valued functions on $\Omega _{n}$
and $\left\Vert f\right\Vert _{E}=\left\Vert f\circ \sigma \right\Vert _{E}$
for every permutation $\sigma $ of $\Omega _{n}.$

The following definitions of $p-$convexity and $q-$concavity can be found in
\cite{LT}:

\noindent For $1\leq p<\infty ,\ E$ is called \emph{p-convex} with constant $%
M<\infty $ if
\[
\Big\Vert\Big(\sum_{i=1}^{n}|v_{i}|^{p}\Big)^{\frac{1}{p}}\Big\Vert_{E}\leq M%
\Big(\sum_{i=1}^{n}\Vert v_{i}\Vert _{E}^{p}\Big)^{\frac{1}{p}}
\]

\noindent for every choice of vectors $v_{1},...,v_{n}\in E$.

\noindent For $1\leq q<\infty ,\ E$ is called \emph{q-concave} with constant
$M<\infty $ if
\[
\Big(\sum_{i=1}^{n}\Vert v_{i}\Vert _{E}^{q}\Big)^{\frac{1}{q}}\leq M\ %
\Big\Vert\Big(\sum_{i=1}^{n}|v_{i}|^{q}\Big)^{\frac{1}{q}}\Big\Vert_{E}
\]

\noindent for every choice of vectors $v_{1},...,v_{n}\in E$.

\medskip The definition of UMD spaces can be found e.g in \cite{B}. We shall
only use the following consequence obtained in \cite{B}: if a Banach space $%
Y $ is UMD then $L_{\psi }\mathcal{\otimes }Id_{Y}$ is a bounded operator on
$L^{2}(\mathbb{R}/2\pi \mathbb{Z},Y)$, where $L_{\psi }$ is the Hilbert
transform, i.e. the convolution operator on $L^{2}(\mathbb{R}/2\pi \mathbb{Z)%
}$ defined by $\psi (\theta )=p.v.\cot \frac{\theta }{2}.$ (Actually this
property is equivalent to UMD).

\noindent We shall denote by $K_{E}$ the non commutative Khintchine upper
constant (as recalled in the proof of lemma \ref{lem:ineqIsoGen}).

\begin{myTheorem}
\label{theorem:main} Let $E$ be a symmetric function space on $\Omega _{n}$.
We assume that

either (i) $E$ is 2-concave with constant $1$.

or (ii) $E$ is 2-convex with constant $1$ and q-concave with constant $1$
for some $q\geq 2.$

We accordingly denote

$K_{E}=1$ in case (i)

$K_{E}=Cq$ in case (ii), where $C$ is a universal constant. If $E=L^{q},q>2,$
one may even take $K_{E}=C\sqrt{q}.$

Let $\varphi \in L^{1}([0,\frac{\pi }{2}],d\theta )$. Then for every
function $f:\Omega _{n}\rightarrow \mathbb{C}$
\begin{equation}
\big\Vert\int_{0}^{\frac{\pi }{2}}\varphi (\theta )\frac{d}{d\theta }\cos ^{%
\frac{\Delta }{4}}\theta (f)d\theta \big\Vert_{E}\leq \frac{1}{2}K_{E}%
\big\Vert\,\left\vert \nabla f\right\vert \,\big\Vert_{E}\int_{0}^{\frac{\pi
}{2}}\left\vert \varphi (\theta )\right\vert d\theta .
\label{equ:mainInequality}
\end{equation}

In particular

\begin{enumerate}
\item Poincar\'{e} inequality holds in $E$: \
\[
\big\Vert f-\mathbb{E}f\big\Vert_{E}\leq \frac{\pi }{4}K_{E}\big\Vert %
\,\left\vert \nabla f\right\vert \,\big\Vert_{E}.
\]

\item Let $\theta _{0}\in \lbrack 0,\frac{\pi }{2}[,$ then
\[
\big\Vert f-\cos ^{\frac{\Delta }{4}}\theta _{0}(f)\big\Vert_{E}\leq \frac{1%
}{2}\theta _{0}K_{E}\,\big\Vert\,\left\vert \nabla f\right\vert \,\big\Vert%
_{E}.
\]

\item Let $0<\alpha <\frac{1}{2}.$ Then
\[
\left\Vert \Delta ^{\alpha }(f)\right\Vert _{E}\leq K_{\alpha }K_{E}\big\Vert%
\,\left\vert \nabla f\right\vert \,\big\Vert_{E}
\]%
where $K_{\alpha }=\frac{1}{\Gamma (1-\alpha )}\left\Vert (-Log\cos \theta
)^{-\alpha}\right\Vert _{L^{1}([0,\frac{\pi }{2}])}$.

\item Moreover, if $C_{E}$ is UMD, then
\[
\Vert \Delta ^{\frac{1}{2}}(f)\Vert _{E}\leq H_{E}K_{E}\big\Vert\,\left\vert
\nabla f\right\vert \,\big\Vert_{E}.
\]

Here $H_{E}$ denotes the finite constant
\[
\sqrt{2\pi }H_{E}=\left\Vert L_{\varphi }\mathcal{\otimes }%
Id_{C_{E}}\right\Vert _{L^{2}(\mathbb{R}/2\pi \mathbb{Z},C_{E})\rightarrow
L^{2}(\mathbb{R}/2\pi \mathbb{Z},C_{E})}
\]

where $L_{\varphi }$ is the convolution operator defined on $L^{2}(\mathbb{R}%
/2\pi \mathbb{Z})$ by
\[
\varphi (\theta )=p.v.1_{[-\frac{\pi }{2},\frac{\pi }{2}]}sgn\theta
(-Log\cos \theta )^{-\frac{1}{2}}.
\]
\end{enumerate}
\end{myTheorem}

\begin{example}
The Orlicz space $L^{\Phi }(\Omega _{n}),$ where $\Phi
(x)=x^{2}Log(1+x^{2}), $ $x\geq 0,$ satisfies assumption (ii) above. Indeed,
since $\varphi (x)=xLog(1+x)$ is convex it is an Orlicz function; hence $%
L^{\Phi }(\Omega _{n}),$ being the 2-convexification of $L^{\varphi }(\Omega
_{n}),$ is 2-convex. Since $\psi (x)=x^{\frac{1}{3}}Log(1+x^{\frac{1}{3}})$
is concave, the Luxemburg functional $l^{\psi }$ defined by $l^{\psi
}(g)=\Vert \,|g|\,^{\frac{1}{6}}\Vert _{L^{\Phi }}^{6}$, associated to $\psi
$ satisfies the reverse triangle inequality:%
\[
\Big\Vert(\sum \left\vert f_{i}\right\vert ^{6})^{\frac{1}{6}}\Big\Vert%
_{L^{\Phi }}=(l^{\psi }(\sum \left\vert f_{i}\right\vert ^{6}))^{\frac{1}{6}%
}\geq (\sum l^{\psi }(\left\vert f_{i}\right\vert ^{6}))^{\frac{1}{6}}=(\sum %
\big\Vert f_{i}\big\Vert_{L^{\Phi }}^{6})^{\frac{1}{6}}
\]%
\noindent i.e. $L^{\Phi }(\Omega _{n})$ is 6-concave with constant 1. Then
case $1$ in Theorem \ref{theorem:main} gives an inequality $\left\Vert f-%
\mathbb{E}f\right\Vert _{L^{\Phi }(\Omega _{n})}\leq K\left\Vert
\,\left\vert \nabla f\right\vert \,\right\Vert _{L^{\Phi }(\Omega _{n})},$
as well as $\left\Vert f-\mathbb{E}f\right\Vert _{L^{2}(\Omega _{n})}\leq
K\left\Vert \,\left\vert \nabla f\right\vert \,\right\Vert _{L^{2}(\Omega
_{n})}.$ Both are weaker than Log-Sobolev inequality in the form \cite[%
inequality (4.4)]{BG}:
\[
\big\Vert f-\mathbb{E}f\big\Vert_{L^{\Phi }(\Omega _{n})}\leq \big\Vert%
\,\left\vert \nabla f\right\vert \,\big\Vert_{L^{2}(\Omega _{n})}.
\]
\end{example}

\paragraph*{Organization of the paper} Theorem \ref{theorem:main} will be
proved in section 3, the concentration phenomenon on the cube in
section 4, reverse inequalities (Theorem \ref{5.1}) in section 5.
Section 6 is devoted to similar results for operators in the CAR
algebra. A remark related to theorem \ref{theorem:main} is given in
the Appendix.

\paragraph*{Acknowledgements}
The authors are grateful to Assaf Naor and Gideon Schechtman for
sub-section $5.5$, and for giving the permission to add it to the
paper. The first author wants to thank Gideon Schechtman for helpful
discussions, and Keith Ball for introducing her to the subject and
for helpful discussions.

\section{\textbf{Preliminaries}}

For the proof of theorem \ref{theorem:main} we need more notation and some
facts, which we take from \cite{LP}.

\subsection{\textbf{\ A *- representation of }$L^{\infty }(\Omega _{n})$
\textbf{into }$\mathcal{M}_{2^{n}}$}

The *-algebra $\mathcal{M}_{2}$ of $2\times 2$ matrices is linearly spanned
by $Id$ and the hermitian unitary Pauli matrices

\[
U=\left(
\begin{array}{cc}
1 & 0 \\
0 & -1%
\end{array}
\right) ,Q=\left(
\begin{array}{cc}
0 & 1 \\
1 & 0%
\end{array}
\right) ,P=\left(
\begin{array}{cc}
0 & i \\
-i & 0%
\end{array}
\right) .
\]

\noindent Note that

\[
QP=-PQ=-iU,
\]

\noindent so that $\mathcal{M}_{2}$ is the *-algebra spanned by $P,Q.$ The
*-algebra $M_{2}$ spanned by $Q$ is *-isomorphic to $L^{\infty }(\Omega
_{1}) $ since $Q$ is the matrix of the pointwise multiplier $\omega $ acting
on $L^{2}(\Omega _{1})$ equipped with the o.n basis $1,\omega .$

The *-algebra $\mathcal{M}_{2^{n}}$ of $2^{n}\times 2^{n}$ matrices equipped
with its unique normalized trace $\tau _{n}$ is identified with $\mathcal{M}%
_{2}\otimes ..\otimes \mathcal{M}_{2}$ ($n$ times) equipped with $\tau
_{1}\otimes ..\otimes \tau _{1}.$ Let

\[
Q_{j}=Id\otimes ..\otimes Q\otimes Id\otimes ..\otimes Id,\;
\]

\noindent where $Q$ occurs in the j$^{th}$ factor. Let $M_{n}$ be the
*-subalgebra spanned by the $Q_{j}$'s, $1\leq j\leq n.$ We denote the
canonical embedding by

\[
J_{n}:M_{n}\rightarrow \mathcal{M}_{2^{n}}.
\]

\noindent Let

\[
P_{j}=Id\otimes ..\otimes P\otimes Id\otimes ..\otimes Id,\;
\]

\noindent so that for $j\neq k$

\[
Q_{j}Q_{k}=Q_{k}Q_{j},\ P_{j}P_{k}=P_{k}P_{j},\ Q_{j}P_{k}=P_{k}Q_{j}
\]
\[
Q_{j}P_{j}=-P_{j}Q_{j}.
\]

\noindent For $A,B\subset \{1,..,n\},\ A=\{i_{1},..,i_{k}\},\
B=\{j_{1},..,j_{m}\},\ i_{1}<..<i_{k},\ j_{1}<..<j_{m},$ let

\[
P_{B}Q_A = P_{j_{1}}...P_{j_{m}} Q_{i_{1}}...Q_{i_{k}}
\]

\noindent and $P_{\varnothing }=Q_{\varnothing }=Id$.

\noindent The matrices $\{P_B Q_{A}\}$ linearly span $\mathcal{M}_{2^{n}}$.

\begin{myLemma}
\label{lemma:traceZero} Let $A,\ B\subset \{1,..,n\},\ A$ or $B\neq
\emptyset .$ Then
\[
\tau _{n}(P_{B}Q_{A})=0.
\]
\end{myLemma}

\begin{proof}
The matrix $P_{B}Q_{A}$ is a tensor product of $2\times 2$ matrices $%
C_{1}\otimes ...\otimes C_{n}.$ At least one of these belongs to $%
\{Q,P,iU\}. $ \noindent \noindent Then $\tau
_{n}(P_{B}Q_{A})=\prod_{j=1}^{n}\tau _{1}(C_{j})=0$.
\end{proof}

Note that $Q_{j}$ is the matrix of the multiplier $\omega _{j}$ acting on $%
L^{2}(\Omega _{n}),$ when $L^{2}(\Omega _{n})$ is equipped with the o.n
Walsh basis$,$ hence $Q_{A}$ is the matrix of the multiplier $\omega _{A}.$
We denote by $\mathcal{I}_{n}:L^{\infty }(\Omega _{n})\rightarrow \mathcal{M}
_{2^{n}}$ the *-representation such that $\omega _{j}\rightarrow Q_{j},$
hence $M_{n}=\mathcal{I}_{n}(L^{\infty }(\Omega _{n})).$ $\mathcal{I} _{n}$
induces operators on $M_{n}$ corresponding to $D_{j}$ and $N$ and we still
denote them by the same letters$,$ in particular the annihilation operator

\[
D_{j}(Q_{A})=\left\{
\begin{array}{ll}
0 & \mbox{if $j \not\in A$} \\
Q_{j}Q_{A} & \mbox{if $j \in A$}%
\end{array}
\right. ,
\]

\noindent the creation operator

\[
D^{\ast}_j(Q_A) = \left\{
\begin{array}{ll}
Q_j Q_A & \mbox{if $j \not\in A$} \\
0 & \mbox{if $j \in A$}%
\end{array}
\right. ,
\]

\noindent and the number operator

\[
N(Q_{A})=\sum_i D^{\ast}_i D_i(Q_A) = \left\vert A\right\vert Q_{A},\ \
N(Id)=0.
\]

\noindent So, with our abuse of notation,

\[
D_{j}(\mathcal{I}_{n}(f))=\mathcal{I}_{n}(D_{j}(f)),\;f\in L^{\infty
}(\Omega_{n}).
\]

\bigskip

\subsection{\textbf{\ A one parameter group of inner *-automorphisms of }$%
\mathcal{M}_{2^{n}}$}

\noindent Let

\[
\mathcal{R}_{\theta }=R_{\theta }\otimes ..\otimes R_{\theta }\in \mathcal{M}
_{2^{n}},\ \ R_{\theta }=\left(
\begin{array}{cc}
1 & 0 \\
0 & e^{i\theta }%
\end{array}
\right) ,\ \ \theta \in \mathbb{R}.
\]

\noindent Since $\mathcal{R}_{\theta }$ is unitary, the action $T\rightarrow
\mathcal{R}_{\theta }^{\ast }T\mathcal{R}_{\theta }$ is an inner
*-automorphism of $\mathcal{M}_{2^{n}}$ which preserves the trace $\tau _{n}$
. Since $R_{\theta }^{\ast }QR_{\theta }=\cos \theta \,Q+\sin \theta \,P$ in
$\mathcal{M}_{2},$

\[
\mathcal{R}_{\theta }^{\ast }Q_{j}\mathcal{R}_{\theta }=\cos \theta
\,Q_{j}+\sin \theta \,P_{j}.
\]

\noindent In particular

\[
P_{j}=\mathcal{R}_{\frac{\pi }{2}}^{\ast }Q_{j}\mathcal{R}_{\frac{\pi }{2}}
\]

\noindent and

\begin{equation}
\mathcal{R}_{\theta }^{\ast }Q_{A}\mathcal{R}_{\theta }=\prod\limits_{j\in
A}(\cos \theta \,Q_{j}+\sin \theta \,P_{j}).  \label{equ:rTheta}
\end{equation}

\noindent Clearly the set of these automorphisms of $\mathcal{M}_{2^{n}}$
(when $\theta $ runs through $\mathbb{R}$ or $\mathbb{R}/2\pi \mathbb{Z})$
is a one parameter group; let $\mathcal{D}$ be its generator. In other
words, we may denote

\[
e^{\theta \mathcal{D}}(T) = \mathcal{R}_{\theta }^{\ast }T\mathcal{R}
_{\theta} ,\ T\in \mathcal{M}_{2^{n}}.
\]

\noindent Observation: $\mathcal{R}_{\theta }$ can be written as $e^{\theta
A},\theta \in \mathbb{R},$ for some (diagonal) antisymmetric matrix $A\in
\mathcal{M}_{2^{n}}$. Then $\mathcal{D}$ is the *-inner derivation defined
by $A,$ i.e. \ $\mathcal{D}(T)=-[A,T]=-AT+TA$. Indeed,

\[
\lim_{\theta \rightarrow 0}\frac{e^{-\theta A}Te^{\theta A}-T}{\theta }
=-AT+TA.
\]

\subsection{\textbf{The Schatten space }$C_{E}$}

Let $E$ be a symmetric function space on $\Omega _{n}$ as defined in the
preliminaries. Let $C_{E}(\mathcal{M}_{2^{n}},\tau _{n})$ be the Schatten
space associated to $E$ and the Von Neumann algebra $\mathcal{M}_{2^{n}}:$
for $S\in \mathcal{M}_{2^{n}},$

\[
\left\Vert S\right\Vert _{C_E}=\left\Vert (s_{k}(S))_{1\leq k\leq
2^{n}}\right\Vert _{E},
\]

\noindent where $(s_{k}(S))_{1\leq k\leq 2^{n}}$ is the sequence of
eigenvalues of $\left\vert S\right\vert=(S^{\ast}S)^{\frac{1}{2}},$ written
in decreasing order. When $E=L^{p}(\Omega _{n}), 1 \leq p < \infty,\ C_{E}$
is denoted by $C_{p}$.

\noindent The $C_{E}$ norm is unitarily invariant. In particular $e^{\theta
\mathcal{D}}$ is an isometry for all Schatten norms $C_{E}.$

We claim that the *- representation $\mathcal{I}_{n}$ of $L^{\infty }(\Omega
_{n})$ into $\mathcal{M}_{2^{n}}$ (such that $\omega _{j}\rightarrow Q_{j})$
extends as an isometry: $E\rightarrow C_{E}(\mathcal{M}_{2^{n}},\tau _{n}).$
Indeed $\mathcal{I}_{n}=\mathcal{I}_{1}\otimes ..\otimes \mathcal{I}_{1}$ and

\[
\mathcal{I}_{1}(\omega )=Q=\rho ^{\ast }U\rho ,\;\rho =\frac{1}{\sqrt{2}}
\left(
\begin{array}{cc}
1 & 1 \\
-1 & 1%
\end{array}
\right) .
\]

\noindent Since $\rho $ is unitary and $U$ is diagonal, the * -automorphism $%
\mathcal{V}_{n}$ of $\mathcal{M}_{2^{n}}$ defined by
\[
\mathcal{V}_{n}(T)=(\rho \otimes ..\otimes \rho )T(\rho ^{\ast }\otimes
..\otimes \rho ^{\ast })
\]

\noindent sends $M_{n}$ onto the *-algebra $\mathcal{A}_{n}$ of diagonal $%
2^{n}\times 2^{n}$ matrices, which can be identified with the algebra of
functions on $\{1,..,2^{n}\}.$ In particular, for $f\in E,\ \mathcal{I}
_{n}(f)$ and $\mathcal{V}_{n}\mathcal{I}_{n}(f)$ have the same $C_{E}$ norm.
On the other hand, the *-representation $\mathcal{V}_{n}\mathcal{I}
_{n}:L^{\infty }(\Omega _{n})\rightarrow \mathcal{A}_{n}$ has the form $%
f\rightarrow f\circ \varphi $ where $\varphi $ is a bijection from $\Omega
_{n}$ onto $\{1,..,2^{n}\}.$ By definition, the $C_{E}$ norm of the diagonal
matrix whose diagonal entries are $f\circ \varphi $ is $\left\Vert
f\right\Vert _{E}.$

\subsection{\textbf{The conditional expectation} $\mathcal{E}_{M_{n}}$}

\noindent Let $\mathcal{E}_{M_{n}}$ be the conditional expectation: $%
\mathcal{M}_{2^{n}}\rightarrow M_{n}$. For $S\in \mathcal{M}_{2^{n}}$, $%
\mathcal{E}_{M_{n}}(S)$ is the unique element of $M_{n}$ such that $\tau
_{n}(V\mathcal{E}_{M_{n}}(S))=\tau _{n}(VS)$ for all $V\in M_{n}$,
equivalently $\tau _{n}(Q_{A}\mathcal{E}_{M_{n}}(S))=\tau _{n}(Q_{A}S)$ for
every $A\subset \{1,..,n\}.$

\medskip

The following is the factorization formula of \cite[lemma 3.1 a)]{LP}.

\begin{myLemma}
\label{lem:factFormula}. For all $T\in M_{n}$
\begin{equation}
\cos ^{N}\theta (T)=\mathcal{E}_{M_{n}}e^{\theta \mathcal{D}}J_{n}(T)=\frac{1%
}{2}\mathcal{E}_{M_{n}}(e^{\theta \mathcal{D}}+e^{-\theta \mathcal{D}%
})J_{n}(T).  \label{equ:factorFormula}
\end{equation}
\end{myLemma}

\begin{proof}
We observe that for $A\subset \{1,..,n\}$ and non empty $B\subseteq
\{1,..,n\}$
\[
\mathcal{E}_{M_{n}}(P_{B}Q_{A})=0.
\]%
Indeed, for $A^{\prime }\subseteq \{1,..,n\},\ \tau _{n}(Q_{A^{\prime
}}P_{B}Q_{A})=0$ by lemma \ref{lemma:traceZero}. Thus expanding the product
\[
e^{\theta \mathcal{D}}=\prod\limits_{j\in A}(\cos \theta \,Q_{j}+\sin \theta
\,P_{j})
\]%
and using the commutation relations, one gets the result.
\end{proof}

Equation (\ref{equ:factorFormula}) expresses a unitary dilation of the
contraction $\cos ^{N}\theta $ acting on the space $C_{2}=L^{2}(M_{n},\tau
_{n})$ (the Hilbert Schmidt operators on $L^{2}(\Omega _{n}))$ or, as well,
expresses the contraction $\cos ^{N}\theta $ acting on $M_{n}$ as the
compression of a *- automorphism of $\mathcal{M}_{2^{n}}$ (Springsteen
factorization). It is an analogue of Mehler formula for $\cos ^{L}\theta $ \
where $L$ is the Ornstein-Uhlenbeck operator in the gaussian setting.

\begin{myLemma}
\label{lem:condException} Let $E$ be a symmetric function space on $\Omega
_{n}.$ The conditional expectation $\mathcal{E}_{M_{n}}$ is a contraction of
$C_{E}(\mathcal{M}_{2^{n}})$.
\end{myLemma}

\begin{proof}
It suffices to verify that $\mathcal{V}_{n}\mathcal{E}_{M_{n}}\mathcal{V}
_{n}^{-1}:\mathcal{M}_{2^{n}}\rightarrow \mathcal{M}_{2^{n}}$ is a
contraction of $C_{E}$. Notice that $\mathcal{V}_{n}\mathcal{E}_{M_{n}}
\mathcal{V}_{n}^{-1}$ is just the restriction of a $2^{n}\times 2^{n}$
matrix to its diagonal part. Let%
\[
U_{j}=Id\otimes ...\otimes U\otimes Id\otimes ..\otimes Id,
\]
\noindent $U$ occurring in the $j^{th}$ factor. As well known, since
\[
U_{j}P_{j}U_{j}=-P_{j},\ U_{j}P_{k}U_{j}=P_{k},\ \ \ k\neq j,
\]%
\[
U_{j}Q_{j}U_{j}=-Q_{j},\ U_{j}Q_{k}U_{j}=Q_{k},\ \ \ k\neq j,
\]%
this restriction can be written as $\mathcal{H}_{n}...\mathcal{H}_{1}$ where
\[
\mathcal{H}_{j}(T)=\frac{1}{2}(T+U_{j}TU_{j}),\;T\in \mathcal{M}_{2^{n}}.
\]
\noindent Every $\mathcal{H}_{j}$ is obviously a contraction of $C_{E}$.

\noindent We shall give another more abstract proof in lemma \ref%
{lem:condExceptionMprime}.

\end{proof}

\section{\textbf{Proof of the main theorem}}

\noindent We shall need the following claims.

\begin{myLemma}
\label{lem:DExpression} For $T\in M_{n},$
\[
\mathcal{D}(T)=\sum_{j=1}^{n}P_{j}D_{j}(T).
\]
\end{myLemma}

\begin{proof}
It is enough to verify the formula for $T=Q_{A}$ and $T=Id.$ It is obvious
in the second case. By definition $\mathcal{D}(Q_{A})$ is the derivative at $%
\theta =0$ of $\theta \rightarrow e^{\theta \mathcal{D}}(Q_{A})=\prod%
\limits_{j\in A}(\cos \theta \,Q_{j}+\sin \theta \,P_{j}).$ Taking into
account the commutation relations, this is
\[
\sum_{j\in A}P_{j}Q_{A\backslash \{j\}}=\sum_{j=1}^{n}P_{j}D_{j}(Q_{A}).
\]
\end{proof}

\begin{myLemma}
\label{lem:DInvSigns} \noindent For every choice of signs $(\varepsilon
_{j})_{1\leq j\leq n}$, there is an inner * -automorphism of $\mathcal{M}
_{2^{n}}$ which maps $\sum_{j=1}^{n}P_{j}D_{j}(T)$ to $\sum_{j=1}^{n}
\varepsilon _{j}P_{j}D_{j}(T),$ $T\in M_{n}.$ In particular, for every
symmetric function space $E$ on $\Omega _{n},$
\[
\Big\Vert\sum_{j=1}^{n}P_{j}D_{j}(T)\Big\Vert_{C_{E}}=\Big\Vert
\sum_{j=1}^{n}\varepsilon _{j}P_{j}D_{j}(T)\Big\Vert_{C_{E}}.
\]
\end{myLemma}

\begin{proof}
The inner * -automorphism of $\mathcal{M}_{2^{n}}$ defined by:
\[
S\rightarrow Q_{j}SQ_{j}
\]
\noindent leaves $M_{n}$ and the $P_{k}$'s ($k\neq j)$ invariant and maps $%
P_{j}$ to $-P_{j}$. For every choice of signs $(\varepsilon _{j})_{1\leq
j\leq n},$ let $A_{\varepsilon }=\{i\mid \varepsilon _{i}=-1,1\leq i\leq
n\}. $ Then the inner * -automorphism:

\[
S\rightarrow Q_{A_{\varepsilon }}SQ_{A_{\varepsilon }}
\]
maps $\sum_{j=1}^{n}P_{j}D_{j}(T)$ to $\sum_{j=1}^{n}\varepsilon
_{j}P_{j}D_{j}(T),$ $T\in M_{n}.$ Since the $C_{E}$ norm is unitarily
invariant, this implies the norm equality.
\end{proof}

\begin{myLemma}
\label{lem:ineqIsoGen} For $T\in M_{n}$ and $E,K_{E}$ as in theorem \textbf{%
\ref{theorem:main}}
\[
\Big\Vert\mathcal{D}(T)\Big\Vert_{C_{E}}\leq K_{E}\Big\Vert%
(\sum_{j=1}^{n}\left\vert D_{j}(T)\right\vert ^{2})^{\frac{1}{2}}\Big\Vert%
_{C_{E}}.
\]
\end{myLemma}

\begin{proof}
By lemmas \ref{lem:DExpression} and \ref{lem:DInvSigns},
\[
\Big\Vert \mathcal{D}(T)\Big\Vert _{C_{E}}=\mathbb{E}\Big\Vert
\sum_{j=1}^{n}\varepsilon _{j}P_{j}D_{j}(T)\Big\Vert_{C_{E}}
\]

\noindent where $(\varepsilon _{j})_{1\leq j\leq N}$ denote iid Bernoulli
variables.

We consider two cases:

\paragraph*{Case 1:}
$E$ is $2$-convex and $q$-concave with constant 1.

\noindent By the non-commutative Khintchine upper inequality we have
\[
\mathbb{E}\Big\Vert\sum_{j=1}^{n}\varepsilon _{j}P_{j}D_{j}(T)\Big\Vert %
_{C_{E}}
\]
\[
\leq K_{E}\max \Big\{\Big\Vert(\sum_{j=1}^{n}P_{j}D_{j}(T)D_{j}(T)^{\ast
}P_{j})^{\frac{1}{2}}\Big\Vert_{C_{E}},\Big\Vert(\sum_{j=1}^{n}D_{j}(T)^{
\ast }D_{j}(T))^{\frac{1}{2}}\Big\Vert_{C_{E}}\Big\}.
\]

\noindent For $C_{E}=C_{p},2\leq p<\infty ,$ see \cite{LP'} and \cite[p. 106]%
{P'} for a better constant; for the general 2-convex and q-concave case, see
\cite[theorem 1.3]{LPX}.

\noindent By the commutation relations
\[
P_{j}D_{j}(T)D_{j}(T)^{\ast }P_{j}=D_{j}(T)D_{j}(T)^{\ast }=D_{j}(T)^{\ast
}D_{j}(T).
\]

\noindent Indeed, for $T=Q_{A},$ $D_{j}(Q_{A})=Q_{j}Q_{A}$ does not involve $%
Q_{j},$ hence commutes with $P_{j}$. \noindent The second equality holds
since $M_{n}$ is a commutative algebra. Hence
\[
\Big\Vert \mathcal{D}(T) \Big\Vert _{C_{E}}\leq K_{E}\Big\Vert%
(\sum_{j=1}^{n}\left\vert D_{j}(T)\right\vert ^{2})^{\frac{1}{2}}\Big\Vert%
_{C_{E}}.
\]

\paragraph*{Case 2:}
E is $2$-concave.

\noindent By the easy part of non-commutative Khintchine inequality (whose
proof is the same as for $C_{p},$ $1\leq p<2$, see \cite{LPP}) we have:
\[
\mathbb{E}\Big\Vert\sum_{j=1}^{n}\varepsilon _{j}P_{j}D_{j}(T)\Big\Vert%
_{C_{E}}\leq \Big\Vert(\sum_{j=1}^{n}D_{j}(T)^{\ast }P_{j}^{\ast
}P_{j}D_{j}(T))^{\frac{1}{2}}\Big\Vert_{C_{E}}=\Big\Vert(\sum_{j=1}^{n}\left%
\vert D_{j}(T)\right\vert ^{2})^{\frac{1}{2}}\Big\Vert_{C_{E}}.
\]
\end{proof}

\begin{myLemma}
\label{lemma:mainInequality} Let $E$ be a symmetric function space on $%
\Omega _{n}$ and let $\varphi \in L^{1}([0,\frac{\pi }{2}]).$ For $T\in
M_{n} $
\[
\Big\Vert\int_{0}^{\frac{\pi }{2}}\varphi (\theta )\frac{d}{d\theta }\cos
^{N}\theta \,(T)d\theta \Big\Vert_{C_{E}}\leq \left\Vert \varphi \right\Vert
_{L^{1}([0,\frac{\pi }{2}])}\Big\Vert\mathcal{D}(T)\Big\Vert_{C_{E}}.
\]

The constant $\left\Vert \varphi \right\Vert _{1}$ may be replaced by $%
\left\Vert L_{\varphi }\otimes Id_{C_{E}}\right\Vert _{L^{2}(\mathbb{R}/2\pi
\mathbb{Z},C_{E})\rightarrow L^{2}(\mathbb{R}/2\pi \mathbb{Z},C_{E})},$
where $L_{\varphi }$ is the convolution operator on $L^{2}(\mathbb{R}/2\pi
\mathbb{Z})$ defined by $1_{[0,\frac{\pi }{2}]}\varphi .$
\end{myLemma}

\begin{proof}
Using factorization (\ref{equ:factorFormula}),
\[
\int_{0}^{\frac{\pi }{2}}\varphi (\theta )\frac{d}{d\theta }\cos ^{N}\theta
\,(T)d\theta =\int_{0}^{\frac{\pi }{2}}\varphi (\theta )\frac{d}{d\theta }%
\mathcal{E}_{M_{n}}e^{\theta \mathcal{D}}(T)d\theta =\mathcal{E}%
_{M_{n}}\int_{0}^{\frac{\pi }{2}}\varphi (\theta )e^{\theta \mathcal{D}}%
\mathcal{D}(T)d\theta .
\]%
\noindent Then, using the contractivity of conditional expectations on
Schatten spaces and the isometry $e^{\theta \mathcal{D}}$,
\begin{eqnarray*}
\Big\Vert\mathcal{E}_{M_{n}}\int_{0}^{\frac{\pi }{2}}\varphi (\theta
)e^{\theta \mathcal{D}}\mathcal{D}(T)d\theta \Big\Vert_{C_{E}} &\leq &%
\Big\Vert\int_{0}^{\frac{\pi }{2}}e^{\theta \mathcal{D}}\mathcal{D}%
(T)\varphi (\theta )d\theta \Big\Vert_{C_{E}} \\
&\leq &\int_{0}^{\frac{\pi }{2}}\Big\Vert\mathcal{\ D}(T)\Big\Vert%
_{C_{E}}\left\vert \varphi (\theta )\right\vert d\theta .
\end{eqnarray*}%
The last assertion of the lemma follows from the transference theorem (see
e.g. \cite[Theorem 2.8]{BGM}): since $\theta \rightarrow e^{\theta \mathcal{D%
}}$ is a representation of the torus $\mathbb{R}/2\pi \mathbb{Z}$ into
isometries of $C_{E},$%
\[
\Big\Vert\int_{0}^{\frac{\pi }{2}}e^{\theta \mathcal{D}}\mathcal{D}%
(T)\varphi (\theta )d\theta \Big\Vert_{C_{E}}\leq \left\Vert L_{\varphi
}\otimes Id_{C_{E}}\right\Vert _{L^{2}(\mathbb{R}/2\pi \mathbb{Z}%
,C_{E})\rightarrow L^{2}(\mathbb{R}/2\pi \mathbb{Z},C_{E})}\Big\Vert\mathcal{%
\ D}(T)\Big\Vert_{C_{E}}.
\]
\end{proof}

\noindent \textbf{Proof of theorem \ref{theorem:main}}: Since $\mathcal{I}%
_{n}$ is an isometry: $E\rightarrow C_{E},$ inequality (\ref%
{equ:mainInequality}) is now rewritten as: for $T\in M_{n}$
\[
\Big\Vert\int_{0}^{\frac{\pi }{2}}\varphi (\theta )\frac{d}{d\theta }\cos
^{N}\theta (T)d\theta \Big\Vert_{C_{E}}\leq K_{E}\int_{0}^{\frac{\pi }{2}%
}\left\vert \varphi (\theta )\right\vert d\theta \Big\Vert%
(\sum_{j=1}^{n}\left\vert D_{j}(T)\right\vert ^{2})^{\frac{1}{2}}\Big\Vert%
_{C_{E}}.
\]

\noindent This follows from lemmas \ref{lem:ineqIsoGen} and \ref%
{lemma:mainInequality}.

\paragraph*{Cases 1 and 2:}

\noindent We now take $\varphi (\theta )=1_{[0,\theta _{0}]}.$ If $0<\theta
_{0}<\frac{\pi }{2}$
\[
-\int_{0}^{\theta _{0}}\frac{d}{d\theta }\cos ^{\frac{\Delta }{4}}\theta
(f)d\theta =f-\cos ^{\frac{\Delta }{4}}\theta _{0}(f)
\]

\noindent hence (\ref{equ:mainInequality}) proves case $2.$ If $\theta _{0}=%
\frac{\pi }{2}$, $\cos ^{\frac{\Delta }{4}}\theta _{0}(f)$ is understood as $%
\lim_{\theta \rightarrow \frac{\pi }{2}}\cos ^{\frac{\Delta }{4}}\theta (f).$
By (\ref{+}) this limit is $\mathbb{E}f$, which proves case $1$.

\paragraph*{Case 3:}

$0<\alpha <\frac{1}{2}$.

\noindent For $\lambda >0$ and $\beta >0$%
\begin{equation}
\Gamma (\beta )\lambda ^{-\beta }=\int_{0}^{\infty }e^{-t\lambda }t^{\beta
-1}dt=\int_{0}^{\frac{\pi }{2}}\cos ^{\lambda -1}\theta (-Log\cos \theta
)^{\beta -1}\sin \theta d\theta .  \label{31}
\end{equation}

\noindent Hence, for every function $f$ on $\Omega _{n},$ taking $\beta
=1-\alpha ,\alpha <1,$
\begin{eqnarray}
\Gamma (1-\alpha )N^{\alpha }(f) &=&\Gamma (1-\alpha )N^{\alpha
-1}N(f)=\int_{0}^{\frac{\pi }{2}}\cos ^{N-1}\theta \;N(f)(-Log\cos \theta
)^{-\alpha }\sin \theta d\theta  \nonumber \\
&=&-\int_{0}^{\frac{\pi }{2}}(-Log\cos \theta )^{-\alpha }\frac{d}{d\theta }%
\cos ^{N}\theta (f)d\theta .  \label{****}
\end{eqnarray}

\noindent Since $-Log\cos \theta \sim \frac{1}{2}\theta ^{2}$ when $\theta
\rightarrow 0,$%
\[
\varphi (\theta )=(-Log\cos \theta )^{-\alpha }\in L^{1}([0,\frac{\pi }{2}%
]),\;\;0<\alpha <\frac{1}{2}.
\]

\noindent By (\ref{equ:mainInequality}) and recalling that $N=\frac{\Delta }{%
4},$ this proves case $3$.

Taking $T=\mathcal{I}_{n}(f)\in M_{n},$ and using factorization (\ref%
{equ:factorFormula}) note that we may rewrite (\ref{****}) as%
\begin{equation}
\Gamma (1-\alpha )N^{\alpha }(T)=\mathcal{E}_{M_{n}}\int_{0}^{\frac{\pi }{2}%
}e^{\theta \mathcal{D}}\mathcal{D}(T)(-Log\cos \theta )^{-\alpha }d\theta .
\label{*}
\end{equation}

\paragraph*{Case 4:}

$\alpha =\frac{1}{2}$.

\noindent This result is proved in \cite[theorem 0.1]{LP} \ for $E=$ $%
L^{q}(\Omega _{n})$ with a worse constant. A similar question is considered
in \cite[Proposition 4]{LP''} with the constant below, but in an abstract
setting. By (\ref{****})
\[
\sqrt{\frac{\pi }{2}}N^{\frac{1}{2}}(f)=-\int_{0}^{\frac{\pi }{2}}(-Log\cos
\theta )^{-\frac{1}{2}}\frac{d}{d\theta }\cos ^{N}\theta (f)d\theta .
\]

\noindent Hence, for $T=\mathcal{I}_{n}(f)\in M_{n},$ by (\ref%
{equ:factorFormula}) we get an analogue of (\ref{*}):%
\begin{eqnarray*}
\sqrt{\frac{\pi }{2}}N^{\frac{1}{2}}(T) &=&-\mathcal{E}_{M_{n}}\int_{0}^{%
\frac{\pi }{2}}(-Log\cos \theta )^{-\frac{1}{2}}\frac{e^{\theta \mathcal{D}%
}-e^{-\theta \mathcal{D}}}{2}\mathcal{D}(T)d\theta \\
&=&-\frac{1}{2}\mathcal{E}_{M_{n}}\;pv\int_{-\frac{\pi }{2}}^{\frac{\pi }{2}%
}e^{\theta \mathcal{D}}\mathcal{D}(T)sgn\theta (-Log\cos \theta )^{-\frac{1}{%
2}}d\theta .
\end{eqnarray*}%
As in the proof of lemma \ref{lemma:mainInequality},
\begin{eqnarray*}
\left\Vert N^{\frac{1}{2}}(T)\right\Vert _{C_{E}} &\leq &\sqrt{\frac{1}{2\pi
}}\left\Vert pv\int_{-\frac{\pi }{2}}^{\frac{\pi }{2}}e^{\theta \mathcal{D}}
\mathcal{D}(T)sgn\theta (-Log\cos \theta )^{-\frac{1}{2}}d\theta \right\Vert
_{C_{E}} \\
&\leq &H_{E}\Big\Vert\mathcal{D}(T)\Big\Vert_{C_{E}}.
\end{eqnarray*}

\noindent Let $\varphi (\theta )=1_{[-\frac{\pi }{2},\frac{\pi }{2}%
]}sgn\theta (-Log\cos \theta )^{-\frac{1}{2}}$ and let
\[
\eta =\varphi -\frac{1}{\sqrt{2}}\cot \frac{\theta }{2}\ \;a.s.\;on\;[-\pi
,\pi].
\]

\noindent Since $\eta $ is continuous on $[-\pi ,\pi ]$, $L_{\eta }\otimes
Id_{C_{E}}$ is bounded on $L^{2}(\mathbb{R}/2\pi \mathbb{Z},C_{E})$; since $%
C_{E}$ is UMD, $L_{\varphi }\otimes Id_{C_{E}}$ is bounded too, hence $H_{E}$
is finite.

\noindent By lemma \ref{lem:ineqIsoGen},
\[
\left\Vert N^{\frac{1}{2}}(T)\right\Vert _{C_{E}}\leq H_{E}K_{E}\Big\Vert%
(\sum_{j=1}^{n}\left\vert D_{j}(T)\right\vert ^{2})^{\frac{1}{2}}\Big\Vert%
_{C_{E}}.
\]

\noindent Since $\mathcal{I}_{n}$ is an isometry: $E\rightarrow C_{E},$ this
proves the claim.

\section{\protect\smallskip \textbf{Applications}}

By a standard argument we can now recover Bobkov-G\"{o}tze inequality \cite[%
corollary 2.4]{BG}, namely assertion 1 in the next corollary:

\begin{myCorollary}
\label{cor:exp} For a function $f:\Omega _{n}\rightarrow \mathbb{C}$

\begin{enumerate}
\item
\[
\mathbb{E}e^{\left\vert f-\mathbb{E}f\right\vert }\leq 2\mathbb{E}e^{\frac{%
\pi ^{2}}{32}\left\vert \nabla f\right\vert ^{2}}\leq 2\,e^{\frac{\pi ^{2}}{%
32}\Vert \,\left\vert \nabla f\right\vert \,\Vert _{\infty }^{2}}.
\]

\item Let $0<\alpha <\frac{1}{2}$. Then

\[
\mathbb{E}e^{\left\vert \triangle ^{\alpha }f\right\vert }\leq 2\mathbb{E}e^{%
\frac{1}{2}\,K_{\alpha }^{2}\left\vert \nabla f\right\vert ^{2}}\leq 2\,e^{%
\frac{1}{2}\,K_{\alpha }^{2}\Vert \,\left\vert \nabla f\right\vert \,\Vert
_{\infty }^{2}}
\]

where $K_{\alpha }=\frac{1}{\Gamma (1-\alpha )}\left\Vert (-Log\cos \theta
)^{-\alpha}\right\Vert _{L^{1}([0,\frac{\pi }{2}])}$.
\end{enumerate}
\end{myCorollary}

\begin{proof}
\begin{enumerate}
\item Theorem \ref{theorem:main} case $1$ is applied to $E=L^{2k}(\Omega
_{n});$ in this case the best constant is $%
K_{2k}^{2k}=(2k-1)!!=1.3.5...(2k-1)$ \cite{Bu}. Hence if $\mathbb{E}f=0$
\begin{eqnarray}
\frac{1}{2}\,\mathbb{E}e^{\left\vert f\right\vert } &\leq &\mathbb{E}%
ch\left\vert f\right\vert =1+\sum_{k\geq 1}\frac{1}{(2k)!}\mathbb{E}%
\left\vert f\right\vert ^{2k}=1+\sum_{k\geq 1}\frac{1}{(2k)!}\left\Vert
f\right\Vert _{L^{2k}}^{2k}  \label{100} \\
&\leq &1+\sum_{k\geq 1}\frac{(2k-1)!!}{(2k!)}(\frac{\pi }{4})^{2k}\left\Vert
\,\left\vert \nabla f\right\vert \,\right\Vert _{L^{2k}}^{2k}=\mathbb{E}%
(1+\sum_{k\geq 1}\frac{1}{2^{k}k!}(\frac{\pi }{4})^{2k}\,\left\vert \nabla
f\right\vert ^{2k})  \nonumber \\
&=&\mathbb{E}\exp (\frac{\pi ^{2}}{32}\,\left\vert \nabla f\right\vert
\,^{2}).  \nonumber
\end{eqnarray}

\item Theorem \ref{theorem:main} case $3$ is applied; the computation is the
same as in the first part, replacing $f$ by $\Delta ^{\alpha }f$ in (\ref%
{100}).
\end{enumerate}
\end{proof}

\noindent As well known, corollary \ref{cor:exp} implies a concentration
inequality for the uniform probability $\mathbb{P}_{n}$ on $\Omega _{n}.$

\begin{myCorollary}
\label{cor:conc} For a function $f:\Omega _{n}\rightarrow \mathbb{C}$ and $%
t>0$

\begin{equation}
\mathbb{P}\{\left\vert f-\mathbb{E}f\right\vert >t\}\leq 2\,\exp (-\frac{
8t^{2}}{\pi ^{2}\left\Vert \,\left\vert \nabla f\right\vert \,\right\Vert
_{\infty }^{2}}).  \label{equ:concLap}
\end{equation}

Let $0<\alpha <\frac{1}{2}$. Then,
\begin{equation}
\mathbb{P}\{\left\vert \triangle ^{\alpha }f\right\vert >t\}\leq 2\,\exp (-%
\frac{t^{2}}{2K_{\alpha }^{2}\left\Vert \,\left\vert \nabla f\right\vert
\,\right\Vert _{\infty }^{2}}).  \label{equ:concLapAlpha}
\end{equation}
\end{myCorollary}

\begin{proof}
By Tchebychev inequality and corollary \ref{cor:exp} applied to $\lambda
\left\vert f-\mathbb{E}f\right\vert ,$ $\lambda >0,$
\[
\mathbb{P}\{\left\vert f-\mathbb{E}f\right\vert >t\}\leq e^{-\lambda t}%
\mathbb{E}e^{\lambda \left\vert f-\mathbb{E}f\right\vert }\leq
2\,e^{-\lambda t+\frac{\pi ^{2}}{32}\lambda ^{2}\left\Vert \,\left\vert
\nabla f\right\vert \,\right\Vert _{\infty }^{2}}.
\]%
\noindent Minimizing the right hand side with respect to $\lambda $ gives (%
\ref{equ:concLap}). The same argument gives (\ref{equ:concLapAlpha}).
\end{proof}

\section{\textbf{Reverse Inequalities}}

\subsection{The statement}

We now extend the main result of \cite{LP} \ which dealt with $\beta =\frac{1%
}{2}$ and $E=L^{p}(\Omega _{n}),$ $1<p<\infty $; note that we improve the
constant in case $p>2$. Though the proof is similar we write the details. We
recall that $K_{E}$ (the non commutative Khintchine upper constant of $%
C_{E}) $ and $H_{E}$ are defined in theorem \ref{theorem:main}.

\begin{myTheorem}
\label{5.1}Let $E$ be a symmetric function space on $\Omega _{n}.$ Then, for
$\beta >\frac{1}{2}$ and every function $f$ on $\Omega _{n},$

\begin{enumerate}
\item if \ $E$ is 2-convex (in particular if $E=L^{\infty }(\Omega _{n}))$
\[
\left\Vert \left\vert \nabla f\right\vert \right\Vert _{E}\leq k_{\beta
}\left\Vert \Delta ^{\beta }f\right\Vert _{E}
\]

\item if $E$ is $2$-concave and $r$-convex with $r>1$
\[
\inf_{\partial _{j}(f)=g_{j}+h_{j}}\Big\{ \left\Vert (\sum \left\vert
g_{j}\right\vert ^{2})^{\frac{1}{2}}\right\Vert _{E}+\left\Vert (\sum
\left\vert h_{j}\ast \delta _{e_{j}}\right\vert ^{2})^{\frac{1}{2}
}\right\Vert _{E} \Big\} \leq k_{\beta }K_{E^{\ast }}^{2}\left\Vert \Delta
^{\beta }f\right\Vert _{E}
\]
\end{enumerate}

Here $k_{\beta }=\frac{1}{\Gamma (\beta )}\left\Vert (-Log\cos \theta
)^{\beta -1}\right\Vert _{L^{1}([0,\frac{\pi }{2}])}$.

Moreover, if $C_{E}$ is UMD, similar inequalities hold for $\beta =\frac{1}{2%
},$ replacing $k_{\beta }$ by $H_{E}.$
\end{myTheorem}

\noindent The translate $h_{j}\ast \delta _{e_{j}}$ is defined by $h_{j}\ast
\delta _{e_{j}}(x)=h_{j}(xe_{j}).$

\noindent We shall see in lemma \ref{5.5} that, for $E=L^{p}(\Omega _{n}),$ $%
1<p<2,$ the left hand side of the inequality in case 2 cannot be replaced by
$\left\Vert \left\vert \nabla f\right\vert \right\Vert _{E}$ if $\frac{1}{2}
\leq \beta < \frac{1}{p}$.

\subsection{\textbf{Notation }}

We denote by $\Pi _{j}$ the orthogonal projection of $C_{2}(\mathcal{M}%
_{2^{n}})$ onto its subspace $P_{j}C_{2}(M_{n}),$ $1\leq j\leq n.$ Since the
range of $\Pi _{j}$ is a left and right module over $M_{n},$ if $S\in
\mathcal{M}_{2^{n}}$ and $T\in M_{n},$

\begin{equation}
\Pi _{j}(ST)=\Pi _{j}(S)T,\;\Pi _{j}(TS)=T\Pi _{j}(S).  \label{20}
\end{equation}

\noindent The $\Pi _{j}$'s have pairwise orthogonal ranges by lemma \ref%
{lemma:traceZero}. Let

\[
\Pi =\sum_{k=1}^{n}\Pi _{k}.
\]

We denote by $X^{\ast }$ the dual space of a Banach space $X.$ We recall
that $C_{E}^{\ast }=C_{E^{\ast }},$ the duality being defined by

\[
\tau (T^{\ast }S),\;S\in C_{E},T\in C_{E}^{\ast }.
\]

\subsection{Some lemmas}

Besides some previous results, the proof of Theorem \ref{5.1} will need the
following lemmas.

\begin{myLemma}
\label{5.2}Let $T\in M_{n}$ such that $\tau (T)=0.$ Then, for $\beta >0,$%
\begin{equation}
\Gamma (\beta )D_{j}N^{-\beta }(T)=P_{j}\Pi _{j}\int_{0}^{\frac{\pi }{2}%
}e^{\theta \mathcal{D}}(T)(-Log\cos \theta )^{\beta -1}d\theta ,\;1\leq
j\leq n,  \label{30}
\end{equation}

whence%
\begin{equation}
\Gamma (\beta )\sum P_{j}D_{j}N^{-\beta }(T)=\Pi \int_{0}^{\frac{\pi }{2}%
}e^{\theta \mathcal{D}}(T)(-Log\cos \theta )^{\beta -1}d\theta .  \label{**}
\end{equation}

Moreover $(-Log\cos \theta )^{\beta -1}\in L^{1}([0,\frac{\pi }{2}])$ if $%
\beta >\frac{1}{2}.$
\end{myLemma}

\noindent Formula (\ref{**}) is parallel to (\ref{*}).

\begin{proof}
Formula (\ref{**}) is obvious from (\ref{30}). It is enough to prove (\ref%
{30}) for $T=Q_{A},$ $A\neq \emptyset .$ The left hand side is zero if $%
j\notin A.$ By (3),

\[
P_{j}\Pi _{j}e^{\theta \mathcal{D}}(Q_{A})=P_{j}\Pi _{j}\prod\limits_{k\in
A}(\cos \theta \,Q_{k}+\sin \theta \,P_{k})
\]

\noindent and this is zero if $j\notin A.$ If $j\in A,$ in virtue of the
commutation relations,

\[
P_{j}\Pi _{j}e^{\theta \mathcal{D}}(Q_{A})=P_{j}(\prod\limits_{k\in
A,k<j}\cos \theta \,Q_{k})\sin \theta \,P_{j}(\prod\limits_{k\in A,k>j}\cos
\theta \,Q_{k})
\]%
\[
=\sin \theta \cos ^{\left\vert A\right\vert -1}\theta D_{j}(Q_{A}).
\]

\noindent By (\ref{31})

\[
\Gamma (\beta )D_{j}N^{-\beta }(Q_{A})=\int_{0}^{\frac{\pi
}{2}}(-Log\cos \theta )^{\beta -1}\sin \theta \cos ^{\left\vert
A\right\vert -1} \theta d\theta \;D_{j}(Q_{A}),
\]

\noindent so (\ref{30}) is proved if $j\in A.$
\end{proof}

\begin{myLemma}
\label{5.3}Let $E$ be a 2-convex symmetric sequence space on $\Omega _{n}.$
Let $S\in \mathcal{M}_{2^{n}}$ and $S_{j}=P_{j}\Pi _{j}(S)\in M_{n},$ $1\leq
j\leq n,$ i.e. $\Pi (S)=\sum P_{j}S_{j}.$ Then

\begin{enumerate}
\item $\Big\Vert (\sum \left\vert S_{j}\right\vert ^{2})^{\frac{1}{2} }%
\Big\Vert _{C_{E}(M_{n})}\leq \left\Vert S\right\Vert
_{C_{E}(\mathcal{M} _{2^{n}})}$.

\item $\Big\Vert (\sum \vert S_{j}^{\ast }P_{j}\vert ^{2})^{\frac{1 }{2}}%
\Big\Vert _{C_{E}(M_{n})}\leq \left\Vert S\right\Vert _{C_{E}( \mathcal{M}%
_{2^{n}})}.$
\end{enumerate}
\end{myLemma}

\noindent For $E=L^{p}(\Omega _{n}),p>2,$ this was proved by interpolation
in \cite{LP}.

\noindent Note that $\vert S_{j}^{\ast }P_{j} \vert \in M_{n}$ since $\vert
S_{j}^{\ast }P_{j} \vert ^{2}=P_{j}S_{j}S_{j}^{\ast }P_{j}\in M_{n}$.

\begin{proof}
\begin{enumerate}
\item Let $E_{(2)}$ be the 2-concavification of $E$ (see \cite{LT}). By
assumption this is a Banach space. We denote by $F=E_{(2)}^{\ast }$ its dual
space. Then, considering $F^{(2)},$ the 2-concavification of $F,$ and using (%
\ref{20}), we get

\begin{eqnarray*}
\left\Vert (\sum \left\vert S_{j}\right\vert ^{2})^{\frac{1}{2}}\right\Vert
_{C_{E}(M_{n})}^{2} &=&\left\Vert \sum \left\vert S_{j}\right\vert
^{2}\right\Vert _{C_{E_{(2)}}(M_{n})}=\sup_{\left\Vert TT^{\ast }\right\Vert
_{C_{F}(M_{n})}=1}\tau (TT^{\ast }\sum S_{j}^{\ast }S_{j}) \\
&=&\sup_{\left\Vert T^{\ast }\right\Vert _{C_{F^{(2)}}(M_{n})}=1}\sum
\left\Vert S_{j}T\right\Vert _{C_{2}(M_{n})}^{2}=\sup_{T^{\ast }}\sum
\left\Vert P_{j}S_{j}T\right\Vert _{C_{2}(\mathcal{M}_{2^{n}})}^{2} \\
&=&\sup_{T^{\ast }}\left\Vert \sum P_{j}S_{j}T\right\Vert
_{C_{2}}^{2}=\sup_{T^{\ast }}\left\Vert \Pi (S)T\right\Vert
_{C_{2}}^{2}=\sup_{T^{\ast }}\left\Vert \Pi (ST)\right\Vert _{C_{2}}^{2} \\
&\leq &\sup_{T^{\ast }}\left\Vert ST\right\Vert _{C_{2}(\mathcal{M}%
_{2^{n}})}^{2}\leq \sup_{\left\Vert V^{\ast }\right\Vert _{C_{F^{(2)}}(%
\mathcal{M}_{2^{n}})}=1}\left\Vert SV\right\Vert _{C_{2}(\mathcal{M}%
_{2^{n}})}^{2} \\
&=&\left\Vert \left\vert S\right\vert ^{2}\right\Vert _{C_{E_{(2)}}(\mathcal{%
M}_{2^{n}})}=\left\Vert S\right\Vert _{C_{E}(\mathcal{M}_{2^{n}})}^{2}.
\end{eqnarray*}

\item The proof is similar:

\begin{eqnarray*}
\left\Vert \sum \left\vert S_{j}^{\ast }P_{j}\right\vert ^{2}\right\Vert
_{C_{E_{(2)}}(M_{n})} &=&\sup_{\left\Vert T^{\ast }T\right\Vert
_{C_{F}(M_{n})}=1}\tau (T^{\ast }T\sum P_{j}S_{j}S_{j}^{\ast }P_{j}) \\
&=&\sup_{\left\Vert T\right\Vert _{C_{F^{(2)}}(M_{n})}=1}\sum \left\Vert
TP_{j}S_{j}\right\Vert _{C_{2}(\mathcal{M}_{2^{n}})}^{2}=\sup_{T}\left\Vert
T\Pi (S)\right\Vert _{C_{2}(\mathcal{M}_{2^{n}})}^{2} \\
&=&\sup_{T}\left\Vert \Pi (TS)\right\Vert _{C_{2}(\mathcal{M}%
_{2^{n}})}^{2}\leq \sup_{T}\left\Vert TS\right\Vert _{C_{2}(\mathcal{M}%
_{2^{n}})}^{2}\leq \left\Vert S\right\Vert _{C_{E}(\mathcal{M}_{2^{n}})}^{2}.
\end{eqnarray*}
\end{enumerate}
\end{proof}

The next lemma is proved in \cite{LP} for $E=L^{p}(\Omega _{n}),$ $1<p<2$.
It is needed for the proof of theorem \textbf{\ref{5.1} }only in case 2.

\begin{myLemma}
\label{5.4}Let $E$ be a symmetric function space on $\Omega _{n}.$ Then

\begin{enumerate}
\item for any orthogonal projection on $C_{2}(\mathcal{M}_{2^{n}})$
\[
\left\Vert \Pi \right\Vert _{C_{E}\rightarrow C_{E}}=\left\Vert \Pi
\right\Vert _{C_{E^{\ast }}\rightarrow C_{E^{\ast }}}.
\]

\item if $\ E^{\ast }$ is 2-convex and $q$-concave ($q>1)$

(i) the projection $\Pi $ defined above satisfies
\[
\left\Vert \Pi \right\Vert _{C_{E^{\ast }}\rightarrow C_{E^{\ast }}}\leq
K_{E^{\ast }}
\]

(ii) for $T_{j}\in C_{E}(M_{n})$ and decompositions $T_{j}=V_{j}+W_{j}$ in $%
C_{E}(M_{n}),$
\[
\inf_{T_{j}=V_{j}+W_{j}} \Big\{ \left\Vert (\sum \left\vert V_{j}\right\vert
^{2})^{ \frac{1}{2}}\right\Vert _{C_{E(M_{n})}}+\left\Vert (\sum \left\vert
W_{j}^{\ast }P_{j}\right\vert ^{2})^{\frac{1}{2}}\right\Vert _{C_{E}(M_{n})} %
\Big\} \leq K_{E^{\ast }}\left\Vert \sum P_{j}T_{j}\right\Vert _{C_{E}(%
\mathcal{M}_{2^{n}})}.
\]
\end{enumerate}
\noindent The assumption on $E^{\ast }$ is equivalent to the
assumption that $E$ is $2$-concave and $q^{\prime }$-convex
($\frac{1}{q}+\frac{1}{q^{\prime }}=1).$
\end{myLemma}

\paragraph*{Remark} In the last assertion, one cannot simply apply non commutative
Khintchine inequality in $C_{E}(\mathcal{M}_{2^{n}})$ ( which would
only use the 2-concavity of $E,$ see \cite{LPX}) since we need
decompositions of $T_{j}$ in $M_{n}$ and not decompositions of
$P_{j}T_{j}$ in $\mathcal{M}_{2^{n}}$.

\medskip

\begin{proof}
\begin{enumerate}
\item The first assertion is classical for orthogonal projections:

\[
\left\Vert \Pi (S)\right\Vert _{C_{E}} =\sup_{\left\Vert T\right\Vert
_{C_{E^{\ast }}}=1}\left\vert \tau (T^{\ast }\Pi (S))\right\vert
=\sup_{T}\left\vert \tau (\Pi (T^{\ast })\Pi (S))\right\vert
\]
\[
=\sup_{T}\left\vert \tau (\Pi (T^{\ast })S)\right\vert \leq \left\Vert \Pi
\right\Vert _{C_{E^{\ast }}\rightarrow C_{E^{\ast }}}\left\Vert S\right\Vert
_{C_{E}}.
\]

\item (i) Let $S\in \mathcal{M}_{2^{n}}$ and $\Pi (S)=\sum P_{j}S_{j}.$ By
lemma \ref{lem:DInvSigns}, by the non commutative Khintchine upper
inequality in $C_{E^{\ast }}(\mathcal{M}_{2^{n}})$ (see the proof of lemma %
\ref{lem:ineqIsoGen}, case 1 and by lemma \ref{5.3})

\begin{eqnarray}
\left\Vert \Pi (S)\right\Vert _{C_{E^{\ast }}(\mathcal{M}_{2^{n}})}
&=&\left\Vert \sum P_{j}S_{j}\right\Vert _{C_{E^{\ast }}}=E\left\Vert \sum
\varepsilon _{j}P_{j}S_{j}\right\Vert _{C_{E^{\ast }}}  \nonumber \\
&\leq &K_{E^{\ast }}\max \{\left\Vert (\sum \left\vert S_{j}\right\vert
^{2})^{\frac{1}{2}}\right\Vert _{C_{E^{\ast }}(M_{n})},\left\Vert (\sum
\left\vert S_{j}^{\ast }P_{j}\right\vert ^{2})^{\frac{1}{2}}\right\Vert
_{C_{E^{\ast }}(M_{n})}\}  \label{32} \\
&\leq &K_{E^{\ast }}\left\Vert S\right\Vert _{C_{E^{\ast }}(\mathcal{M}%
_{2^{n}})}.  \nonumber
\end{eqnarray}

(ii) We argue by duality. For $S_{k}\in M_{n},$ $1\leq k\leq n,$ we get

\begin{eqnarray*}
K_{E^{\ast }}\left\Vert \sum P_{j}T_{j}\right\Vert _{C_{E}(\mathcal{M}
_{2^{n}})} &\geq &\sup_{\left\Vert \sum P_{k}S_{k}\right\Vert _{C_{E^{\ast
}}}\leq K_{E^{\ast }}} \Big\vert \tau \Big( (\sum P_{k}S_{k})^{\ast }(\sum
P_{j}T_{j}) \Big) \Big\vert \\
&\geq &\sup_{\left\Vert (\sum \left\vert S_{k}\right\vert ^{2})^{\frac{1}{2}
}\right\Vert _{C_{E^{\ast }}}\leq 1,\left\Vert (\sum \left\vert S_{k}^{\ast
}P_{k}\right\vert ^{2})^{\frac{1}{2}}\right\Vert _{C_{E^{\ast }}}\leq 1} %
\Big\vert \tau \Big( (\sum P_{k}S_{k})^{\ast }(\sum P_{j}T_{j}) \Big) %
\Big\vert \\
&=&\sup_{\left\Vert (\sum \left\vert S_{k}\right\vert ^{2})^{\frac{1}{2}
}\right\Vert _{C_{E^{\ast }}}\leq 1,\left\Vert (\sum \left\vert S_{k}^{\ast
}P_{k}\right\vert ^{2})^{\frac{1}{2}}\right\Vert _{C_{E^{\ast }}}\leq
1}\left\vert \tau (\sum S_{j}{}^{\ast }T_{j})\right\vert =\left\Vert
(T_{j})_{j=1}^{n}\right\Vert
\end{eqnarray*}

\noindent where the second inequality uses (\ref{32}). We shall now explicit
$\left\Vert (T_{j})_{j=1}^{n}\right\Vert $ defined by the previous line.

Let $X$ \ be the Banach space $M_{n}\times ..\times M_{n}$ ($n$ times)
equipped with the norm

\[
\left\Vert (V_{j})_{j=1}^{n}\right\Vert _{X}=\left\Vert (\sum \left\vert
V_{j}\right\vert ^{2})^{\frac{1}{2}}\right\Vert _{C_{E}(M_{n})}.
\]

\noindent As well known, for the duality defined by $\tau (\sum S_{j}^{\ast
}V_{j}),$

\[
\left\Vert (S_{j})_{j=1}^{n}\right\Vert _{X^{\ast }}=\left\Vert (\sum
\left\vert S_{j}\right\vert ^{2})^{\frac{1}{2}}\right\Vert _{C_{E^{\ast
}}(M_{n})}.
\]

\noindent $.$

\noindent Similarly let $Y$ be the Banach space $M_{n}\times ..\times M_{n}$
($n$ times) equipped with the norm

\[
\left\Vert (W_{j})_{j=1}^{n}\right\Vert _{Y}=\left\Vert (P_{j}W_{j}^{\ast
}P_{j})_{j=1}^{n}\right\Vert _{X}=\left\Vert (\sum \left\vert W_{j}^{\ast
}P_{j}\right\vert ^{2})^{\frac{1}{2}}\right\Vert _{C_{E}(M_{n})}.
\]

\noindent (Again, recall that $P_{j}W_{j}^{\ast }P_{j}\in M_{n}).$ We claim
that, for the same duality as above, namely $\tau (\sum S_{j}^{\ast }W_{j}),$

\[
\left\Vert (S_{j})_{j=1}^{n}\right\Vert _{Y^{\ast }}=\left\Vert
(P_{j}S_{j}^{\ast }P_{j})_{j=1}^{n}\right\Vert _{X^{\ast }}=\left\Vert (\sum
\left\vert S_{j}^{\ast }P_{j}\right\vert ^{2})^{\frac{1}{2}}\right\Vert
_{C_{E^{\ast }}(M_{n})}.
\]

\noindent Indeed

\begin{eqnarray*}
\left\Vert (S_{j})_{j=1}^{n}\right\Vert _{Y^{\ast }} &=&\sup_{\left\Vert
(W_{j})_{j=1}^{n}\right\Vert _{Y}=1}\left\vert \tau (\sum S_{j}^{\ast
}W_{j})\right\vert =\sup_{\left\Vert (W_{j})_{j=1}^{n}\right\Vert
_{Y}=1}\left\vert \tau (\sum S_{j}W_{j}^{\ast })\right\vert \\
&=&\sup_{\left\Vert (P_{j}W_{j}^{\ast }P_{j})_{j=1}^{n}\right\Vert
_{X}=1}\left\vert \tau (\sum P_{j}S_{j}P_{j}^{2}W_{j}^{\ast
}P_{j})\right\vert \\
&=&\sup_{\left\Vert (R_{j})_{j=1}^{n}\right\Vert _{X}=1}\left\vert \tau
(\sum P_{j}S_{j}P_{j}R_{j})\right\vert =\left\Vert (P_{j}S_{j}^{\ast
}P_{j})_{j=1}^{n}\right\Vert _{X^{\ast }}.
\end{eqnarray*}

\noindent It follows that

\[
\left\Vert (T_{j})_{j=1}^{n}\right\Vert = \left\Vert
(T_{j})_{j=1}^{n}\right\Vert _{X+Y}=\inf_{S_{j}=V_{j}+W_{j}} \Big\{
\Big\Vert (V_{j})_{j=1}^{n} \Big\Vert _{X}+\left\Vert
(W_{j})_{j=1}^{n}\right\Vert _{Y} \Big\}
\]
\[
= \inf_{S_{j}=V_{j}+W_{j}} \Big\{ \Big\Vert (\sum \left\vert
V_{j}\right\vert ^{2})^{\frac{1}{2}} \Big\Vert _{C_{E}}+\Big\Vert
(\sum \left\vert W_{j}^{\ast }P_{j}\right\vert ^{2})^{\frac{1}{2}} \Big\Vert %
_{C_{E}} \Big\}.
\]

\noindent and this ends the proof of (ii).
\end{enumerate}
\end{proof}

\subsection{\noindent \textbf{Proof of Theorem \protect\ref{5.1} }}

We recall that $e^{\theta \mathcal{D}}$ is an isometry of $C_{E}(\mathcal{M}%
_{2^{n}}).$

\paragraph*{Case 1:}

$E$ is 2-convex.

\noindent It is enough to prove that, if $\mathbb{E}f=0,$

\[
\left\Vert \vert \nabla \Delta ^{-\beta }f \vert \right\Vert _{E}\leq \frac{2%
}{4^{\beta }}k_{\beta }\left\Vert f\right\Vert _{E}
\]

\noindent or, equivalently, for $T=\mathcal{I}_{n}(f),$

\[
\Big\Vert (\sum \vert D_{j}N^{-\beta }(T) \vert ^{2})^{\frac{1}{2}}\Big\Vert %
_{C_{E}(M_{n})} \leq k_{\beta }\left\Vert T\right\Vert _{C_{E}(M_{n})}.
\]

\noindent By (\ref{30}) and lemma \ref{5.3} 1.

\[
\Big\Vert (\sum \vert D_{j}N^{-\beta }(T) \vert ^{2})^{\frac{1}{2} }%
\Big\Vert _{C_{E}} \leq \frac{1}{\Gamma (\beta )} \Big\Vert
\int_{0}^{ \frac{\pi }{2}}e^{\theta \mathcal{D}}(T)(-Log\cos \theta )^{\beta
-1}d\theta \Big\Vert _{C_{E}}
\]
\[
\leq \frac{1}{\Gamma (\beta )}\int_{0}^{\frac{\pi }{2}}(-Log\cos \theta
)^{\beta -1}d\theta \;\left\Vert T\right\Vert _{C_{E}}=k_{\beta }\left\Vert
T\right\Vert _{C_{E}}.
\]

\paragraph*{Case 2:}

$E$ is 2-concave and $r$-convex.

\noindent (\ref{**}) and lemma \ref{5.4} 2. (i) imply

\[
\left\Vert \sum P_{j}D_{j}N^{-\beta }(T)\right\Vert _{C_{E}} = \frac{1}{
\Gamma (\beta )}\Big\Vert \Pi \int_{0}^{\frac{\pi }{2}}e^{\theta \mathcal{D}
}(T)(-Log\cos \theta )^{\beta -1}d\theta \Big\Vert _{C_{E}}
\]
\[
\leq k_{\beta }\left\Vert \Pi \right\Vert _{C_{E}\rightarrow
C_{E}}\left\Vert T\right\Vert _{C_{E}}\leq k_{\beta }K_{E^{\ast }}\left\Vert
T\right\Vert _{C_{E}}.
\]

\noindent Hence

\[
\left\Vert \mathcal{D}(T)\right\Vert _{C_{E}}=\left\Vert \sum
P_{j}D_{j}(T)\right\Vert _{C_{E}}\leq k_{\beta }K_{E^{\ast }} \left\Vert
N^{\beta }T \right\Vert _{C_{E}}.
\]

\noindent Since $E^{\ast }$ is 2-convex and $r^{\prime }-$convex ($\frac{1}{%
r }+\frac{1}{r^{\prime }}=1)$ lemma \ref{5.4} 2. (ii) implies

\[
\inf_{D_{j}(T)=V_{j}+W_{j}} \Big\{ \Big\Vert (\sum \left\vert
V_{j}\right\vert ^{2})^{\frac{1}{2}}\Big\Vert
_{C_{E(M_{n})}}+\Big\Vert (\sum \left\vert W_{j}^{\ast }P_{j}\right\vert
^{2})^{\frac{1}{2}} \Big\Vert _{C_{E}(M_{n})} \Big\}
\]
\[
\leq K_{E^{\ast }}\Big\Vert \sum P_{j}D_{j}(T) \Big\Vert
_{C_{E}}\leq k_{\beta }K_{E^{\ast }}^{2}\Big\Vert N^{\beta }T \Big\Vert %
_{C_{E}}
\]

\noindent where $V_{j},W_{j}\in M_{n}.$

The remaining task is to translate this inequality in $E.$ Note that $%
P_{j}Q_{A}P_{j}=Q_{A}$ if $j\notin A$ and $P_{j}Q_{A}P_{j}=-Q_{A}$ if $j\in
A.$ Hence

\[
P_{j}Q_{A}P_{j}=\mathcal{I}_{n}(\varepsilon _{A}\ast \delta _{e_{j}}).
\]

\noindent For every $W\in M_{n},$ there exists $h\in L^{\infty }(\Omega
_{n}) $ such that $W=\mathcal{I}_{n}(h),$ hence $WW^{\ast }=\mathcal{I}
_{n}(\left\vert h\right\vert ^{2})$ and

\[
\left\vert W_{j}^{\ast }P_{j}\right\vert ^{2}=P_{j}WW^{\ast }P_{j}=\mathcal{I%
}_{n}(\left\vert h\right\vert ^{2}\ast \delta _{e_{j}})=\mathcal{I}%
_{n}(\left\vert h\ast \delta _{e_{j}}\right\vert ^{2}).
\]

\noindent Recalling that $\partial _{j}f=2\omega _{j}D_{j}f,$ we thus get

\[
\inf_{\partial _{j}f=g_{j}+h_{j}} \Big\{ \left\Vert (\sum \left\vert
g_{j}\right\vert ^{2})^{\frac{1}{2}}\right\Vert _{E}+\left\Vert (\sum
\left\vert h_{j}\ast \delta _{e_{j}}\right\vert ^{2})^{\frac{1}{2}
}\right\Vert _{E} \Big\} \leq \frac{2}{4^{\beta }}k_{\beta }K_{E^{\ast
}}^{2}\left\Vert \Delta ^{\beta }f\right\Vert _{E},
\]

\noindent which ends the proof of of case 2.

\paragraph*{The UMD case:}

\noindent As in the proof of theorem \ref{theorem:main}, case 4, (\ref{**})
may be replaced by
\[
\sum P_{j}D_{j}N^{-\frac{1}{2}}(T) = \frac{1}{\sqrt{2\pi }}\Pi \int_{0}^{
\frac{\pi }{2}}(e^{\theta \mathcal{D}}-e^{-\theta \mathcal{D}})(T)(-Log\cos
\theta )^{-\frac{1}{2}}d\theta
\]
\[
= \frac{1}{\sqrt{2\pi }}\Pi pv\int_{-\frac{\pi }{2}}^{\frac{\pi }{2}
}e^{\theta \mathcal{D}}(T)sgn\theta (-Log\cos \theta )^{-\frac{1}{2}}d\theta
\]

\noindent and similarly for (\ref{30}). Moreover

\[
\frac{1}{\sqrt{2\pi }} \Big\Vert pv\int_{-\frac{\pi }{2}}^{\frac{\pi }{2}%
}e^{\theta \mathcal{D}}(T)sgn\theta (-Log\cos \theta )^{-\frac{1}{2}}d\theta %
\Big\Vert _{C_{E}} \leq H_{E}\left\Vert T\right\Vert _{C_{E}}.
\]

\noindent The remaining computations come from cases 1 and 2 above.

\subsection{A remark on theorem \textbf{\protect\ref{5.1}}}

\noindent The first author learnt the following result from A. Naor and G.
Schechtman:

\begin{myLemma}
\label{5.5}Let $E=L^{p}(\Omega _{n}),$ and assume that there is a constant $%
C_{\beta }$ such that, for every $n$ and every function $f$ on $\Omega _{n}$
\begin{equation}
\left\Vert\, \left\vert \nabla f\right\vert\, \right\Vert _{p}\leq C_{\beta
} \Vert \Delta ^{\beta }f \Vert _{p}.  \label{22}
\end{equation}

Then $\beta \geq \frac{1}{p}.$
\end{myLemma}

\noindent The fact that (\ref{22}) is false for $\beta =\frac{1}{2}$ was
proved at the end of \cite{LP}. For the sake of completeness we reproduce
Naor and Schechtman\ proof, which, similarly, relies on an idea of D.
Lamberton. The new ingredient is the following lemma, whose proof for $\beta
=\frac{1}{2}$ is different in \cite{LP}:

\begin{myLemma}
Let $E$ be a Banach space, and let $(T_{t})_{t>0}$ be a strongly continuous
semi-group of contractions of $E.$ Let $A$ be its (densely defined) generator%
$.$ Then for $0<\beta <1$ and every $f$ in the domain of $A$
\begin{equation}
\Vert A^{\beta }f\Vert _{E}\leq 4\left\Vert f\right\Vert _{E}^{1-\beta
}\left\Vert Af\right\Vert _{E}^{\beta }.  \label{21}
\end{equation}
\end{myLemma}

\begin{proof}
For $\lambda >0$ and $0<\beta <1$

\[
C_{\beta }\lambda ^{\beta }=\int_{0}^{\infty }\frac{1-e^{-\lambda t}}{%
t^{\beta +1}}dt,
\]

\noindent where

\[
C_{\beta }=\int_{0}^{\infty }\frac{1-e^{-t}}{t^{\beta +1}}dt\geq \frac{1}{2}%
(\int_{0}^{1}\frac{dt}{t^{\beta }}+\int_{1}^{\infty }\frac{dt}{t^{\beta +1}}%
)=\frac{1}{2\beta (1-\beta )}.
\]

\noindent Hence, for $f$ in the domain of $A$ in $E,$

\[
A^{\beta }f=\frac{1}{C_{\beta }}\int_{0}^{\infty }\frac{f-e^{-tA}f}{t^{\beta
+1}}dt.
\]

\noindent For every $M>0$

\[
\int_{0}^{\infty }\frac{f-e^{-tA}f}{t^{\beta +1}}dt=\int_{0}^{M}(%
\int_{0}^{t}A e^{-A}fdu)\frac{dt}{t^{\beta +1}}+\int_{M}^{\infty }\frac{%
f-e^{-tA}f}{t^{\beta +1}}dt.
\]

\noindent Since $e^{-tA}$ is a contraction on $E$

\[
\left\Vert A^{\beta }f\right\Vert _{E}\leq 2\beta (1-\beta ) \left[%
\left\Vert Af\right\Vert _{E}\frac{M^{1-\beta }}{1-\beta }+\left\Vert
f\right\Vert _{E} \frac{2}{\beta }M^{-\beta }\right].
\]

\noindent By choosing $M=2\left\Vert f\right\Vert _{E}\left\Vert
Af\right\Vert _{E}^{-1}$ we get the result.
\end{proof}

\begin{proofof}{Lemma \ref{5.5}}
Let $f_{n}$ \ be the Riesz product

\[
f_{n}=\prod\limits_{j=1}^{n}(1+\omega _{j})=2^{n}1_{\{1,..,1\}}.
\]

\noindent Then $\Delta f_{n}$ and $\left\vert \nabla f_{n}\right\vert $ can
be computed and are not supported on the point $\{1,..,1\},$ which is due to
the fact that the $\partial _{j}$'s are not local operators. As proved by
Lamberton (see \cite{LP}), for $1\leq p<\infty ,$

\[
\left\Vert \Delta f_{n}\right\Vert _{p}\leq 2^{1+\frac{1}{p}}n\left\Vert
f_{n}\right\Vert _{p}
\]

\noindent and

\[
n^{\frac{1}{p}}\left\Vert f_{n}\right\Vert _{p}\leq \left\Vert \left\vert
\nabla f_{n}\right\vert \right\Vert _{p}.
\]

Since for $t>0$ $e^{-t\Delta }$ is Markovian on $L^{\infty }(\Omega _{n})$
and $\mathbb{E}e^{-t\Delta }f=f$ \ for every $f$, $e^{-t\Delta }$ is a
contraction on $L^{\infty }(\Omega _{n})$ and $L^{1}(\Omega _{n})$. By
interpolation, it a contraction on $L^{p}(\Omega _{n}),$ $1\leq p\leq \infty
$ (and on any symmetric function space $E$ on $\Omega _{n}$). Hence (\ref{21}%
) holds for $A=\Delta $ and every $f.$ Thus, if (\ref{22}) \ holds, for
every $n,$

\[
n^{\frac{1}{p}}\left\Vert f_{n}\right\Vert _{p} \leq \Big\Vert
(\sum_{j=1}^{n}\left\vert \partial _{j}f_{n}\right\vert ^{2})^{\frac{1}{2}}%
\Big\Vert_{p} \leq C_{\beta } \Vert \Delta ^{\beta }f_{n} \Vert _{p}
\]
\[
\leq 4C_{\beta }\left\Vert f_{n}\right\Vert _{p}^{1-\beta }\left\Vert \Delta
f_{n}\right\Vert _{p}^{\beta }\leq 4.2^{(1+\frac{1}{p})\beta }C_{\beta
}n^{\beta }\left\Vert f_{n}\right\Vert _{p}.
\]

\noindent This implies $\beta \geq \frac{1}{p}.$
\end{proofof}

\section{\textbf{Inequalities in the CAR *-algebra}}

\subsection{Notation}

We define in $\mathcal{M}_{2^{n}}$

\[
Q_{j}^{\prime }=U\otimes ..\otimes Q\otimes Id\otimes ..\otimes Id,
\]
\[
P_{j}^{\prime }=U\otimes ..\otimes P\otimes Id\otimes ..\otimes Id
\]

\noindent where $U$ occurs in the $j-1$ first factors. $P_{j}^{\prime }$ and
$Q_{j}^{\prime }$ are hermitian and unitary and satisfy the following
canonical anticommutation relations (CAR): for $1\leq j,k\leq n$

\[
Q_{j}^{\prime }P_{k}^{\prime }=-P_{k}^{\prime }Q_{j}^{\prime },
\]
\[
Q_{j}^{\prime }Q_{k}^{\prime }=-Q_{k}^{\prime }Q_{j}^{\prime }, \ \ \ j\neq
k,
\]
\[
P_{j}^{\prime }P_{k}^{\prime }=-P_{k}^{\prime }P_{j}^{\prime }, \ \ \ j\neq
k.
\]

\noindent For $A,B \subset \{1,..,n\},\ A=\{i_{1},..,i_{k}\},\
B=\{j_{1},..,j_{m}\},\ i_{1}<..<i_{k},\ j_{1}<..<j_{m},$ let

\[
P^{\prime}_{B}Q^{\prime}_A = P^{\prime}_{j_{1}}...P^{\prime}_{j_{m}}
Q^{\prime}_{i_{1}}...Q^{\prime}_{i_{k}}
\]

\noindent and $Q^{\prime}_{\varnothing }=P^{\prime}_{\varnothing }=Id$.

\noindent The matrices $\{P^{\prime}_B Q^{\prime}_{A}\}$ linearly span $%
\mathcal{M}_{2^{n}}$.

\medskip

The purpose of this section is to repeat the main results of the previous
sections, replacing functions $f$ on the cube by elements $T$ of the CAR
*-algebra $M_{n}^{^{\prime }}$ spanned by the $n$ anticommuting fermions $%
Q_{1}^{\prime },..,Q_{n}^{\prime }.$ $M_{n}^{\prime }$ is a non commutative
subalgebra of $\mathcal{M}_{2^{n}}$ (we denote by $J_{n}^{\prime }$ the
embedding) and we consider Schatten norms $C_{E}(\mathcal{M}_{2^{n}},\tau
_{n})$ as before. The conditional expectation: $\mathcal{M}%
_{2^{n}}\rightarrow M_{n}^{\prime }$ is denoted by $\mathcal{E}%
_{M_{n}^{\prime }}.$ The annihilation operator $D_{j}^{\prime }:$ $%
C_{2}(M_{n}^{\prime })\rightarrow C_{2}(M_{n}^{\prime })$ is defined by

\[
D_{j}^{\prime }(Q_{A}^{\prime })=\left\{
\begin{array}{ll}
0 & \mbox{if $j \not\in A$} \\
Q_{j}^{\prime }Q_{A}^{\prime } & \mbox{if $j \in A$}%
\end{array}%
\right. .
\]

\noindent and the number operator $N^{\prime }$ by $N^{\prime
}(Q_{A}^{\prime })=\sum_{j=1}^{n}D_{j}^{\prime \ast }D_{j}^{\prime
}(Q_{A}^{\prime })=\left\vert A\right\vert Q_{A}^{\prime },\ \ N^{\prime
}(Id)=0$.

\noindent We denote by $\left\vert \nabla _{s}T\right\vert $ the following
"symmetrized" sum:
\[
|\nabla _{s}T|=\Big
(\sum_{j=1}^{n}\left\vert D_{j}^{\prime }(T)^{\ast }\right\vert
^{2}+\left\vert D_{j}^{\prime }(T)\right\vert ^{2}\Big)^{\frac{1}{2}}.
\]

\noindent We consider as in subsection 2.2 the inner *-automorphism $%
e^{\theta \mathcal{D}}$ of $\mathcal{M}_{2^{n}}:\ T\rightarrow \mathcal{R}%
_{\theta }^{\ast }T\mathcal{R}_{\theta }$. Since $R_{\theta }^{\ast
}UR_{\theta }=U$ in $\mathcal{M}_{2},$ we still have

\[
\mathcal{R}_{\theta }^{\ast }Q_{A}^{\prime }\mathcal{R}_{\theta
}=\prod\limits_{j\in A}(\cos \theta \,Q_{j}^{\prime }+\sin \theta
\,P_{j}^{\prime }).
\]

\subsection{\textbf{The main theorem for CAR algebras.}}

We now state the analogue of theorem \ref{theorem:main}.

\begin{myTheorem}
\label{theorem:mainCAR} Let $E,K_{E},K_{\alpha },H_{E}$ be as in theorem \ref%
{theorem:main} and let $T$ belong to the *-algebra $M_{n}^{\prime }$ spanned
by $n$ fermions. Then
\[
\;\Big\Vert T-\tau _{n}(T)Id\Big\Vert_{C_{E}}\leq \frac{\pi }{2}K_{E}%
\Big\Vert\left\vert \nabla _{s}T\right\vert \Big\Vert_{C_{E}}
\]

\noindent and
\[
\Big\Vert N^{\prime \alpha }(T)\Big\Vert_{C_{E}}\leq K_{\alpha }K_{E}%
\Big\Vert\left\vert \nabla _{s}T\right\vert \Big\Vert_{C_{E}},\;\;0<\alpha <%
\frac{1}{2}.
\]

\noindent Moreover, if $C_{E}$ is UMD,
\[
\Big\Vert N^{\prime \frac{1}{2}}(T)\Big\Vert_{C_{E}}\leq H_{E}K_{E}\Big\Vert%
\left\vert \nabla _{s}T\right\vert \Big\Vert_{C_{E}}.
\]
\end{myTheorem}

The proof of this theorem needs the counterparts in this setting of the
lemmas from parts 2 and 3. Replacing $P_{j},Q_{j},N,J_{n},\mathcal{E}%
_{M_{n}} $ by $Q_{j}^{\prime },P_{j}^{\prime },N^{\prime },J_{n}^{\prime },%
\mathcal{E}_{M_{n}^{\prime }}$ respectively, give with the same proofs the
analogues of lemmas \ref{lemma:traceZero}, \ref{lem:factFormula}, \ref%
{lem:DExpression}. In particular for every $T\in M_{n}^{\prime }$
\[
\mathcal{E}_{M_{n}^{\prime }}e^{\theta \mathcal{D}}J_{n}^{\prime }(T)=\cos
^{N^{\prime }}\theta (T)
\]

\noindent and
\[
\mathcal{D}(T)=\sum_{j=1}^{n}P_{j}^{\prime }D_{j}^{\prime }(T).
\]

\begin{myLemma}
For any $T\in M_{n}^{\prime }$
\[
\lim_{\theta \rightarrow \frac{\pi }{2}}\cos ^{N^{\prime }}\theta (T)=\tau
_{n}(T)Id.
\]
\end{myLemma}

\begin{proof}
We can expand T as follows:
\[
T=\tau _{n}(T)\cdot Id+\sum_{A\subset \{1,..,n\},A\neq \varnothing }\alpha
_{A}Q_{A}^{\prime }.
\]%
\noindent hence
\[
\cos ^{N^{\prime }}\theta (T)=\tau _{n}(T)\cdot Id+\sum_{A\subset
\{1,..,n\},A\neq \varnothing }\alpha _{A}\cos ^{\left\vert A\right\vert
}\theta Q_{A}^{\prime }\rightarrow _{\theta \rightarrow \frac{\pi }{2}}\tau
_{n}(T)Id.
\]
\end{proof}

\begin{myLemma}
\label{lem:condExceptionMprime} For every symmetric function space $E$ on $%
\Omega _{n},$ the conditional expectation $\mathcal{E}_{M_{n}^{\prime }}$ is
a contraction of $C_{E}$.
\end{myLemma}

\begin{proof}
As well known, every conditional expectation is a contraction of $C_{1}$ and
$C_{\infty }.$ By \cite[corollary 2.11]{A} $C_{E}$ is an interpolation space
between $C_{1}$ and $C_{\infty },$ meaning that every operator $A$ which is
bounded both on $C_{1}$ and $C_{\infty }$ is bounded on $C_{E};$ more
precisely
\[
\big\Vert A\big\Vert_{C_{E}\rightarrow C_{E}}\leq \max\{\big\Vert A\big\Vert%
_{C_{1}\rightarrow C_{1}},\big\Vert A\big\Vert_{C_{\infty }\rightarrow
C_{\infty }}\}.
\]
\end{proof}

\begin{myLemma}
\label{lem:DInvSignsCAR} \noindent For every choice of signs $(\varepsilon
_{j})_{1\leq j\leq n}$, there is an inner * -automorphism of $\mathcal{M}%
_{2^{n}}$ which maps $\sum_{j=1}^{n}P_{j}^{\prime }D_{j}^{\prime }(T)$ to $%
\sum_{j=1}^{n}\varepsilon _{j}P_{j}^{\prime }D_{j}^{\prime }(T),\ T\in
M_{n}^{\prime }$ . In particular, for every symmetric function space $E$ on $%
\Omega _{n},$
\[
\Big\Vert\sum_{j=1}^{n}P_{j}^{\prime }D_{j}^{\prime }(T)\Big\Vert_{C_{E}}=%
\Big\Vert\sum_{j=1}^{n}\varepsilon _{j}P_{j}^{\prime }D_{j}^{\prime }(T)%
\Big\Vert_{C_{E}}.
\]
\end{myLemma}

\begin{proof}
It is enough to exhibit for every $j$ an inner * -automorphism of $\mathcal{M%
}_{2^{n}}$ which stabilizes $Q_{i}^{\prime },1\leq i\leq n,$ and $%
P_{k}^{\prime },k\neq j,$ and maps $P_{j}^{\prime }$ to $-P_{j}^{\prime }.$
Let $V_{j}$ be the hermitian and unitary matrix
\[
V_{j}=Id\otimes ..\otimes Q\otimes U\otimes ..\otimes U\in \mathcal{M}%
_{2^{n}}
\]%
\noindent where $Q$ occurs in the $j^{th}$ factor. In virtue of the
commutation relations in $\mathcal{M}_{2},$ in particular $QPQ=-P,$ $QUQ=-U,$
$UQU=-Q,$ $UPU=-P,$
\[
V_{j}Q_{k}^{\prime }V_{j}=Q_{k}^{\prime },\ \ 1\leq k\leq n,
\]%
\[
V_{j}P_{k}^{\prime }V_{j}=P_{k}^{\prime },\ \ k\neq j,\ \ \
V_{j}P_{j}^{\prime }V_{j}=-P_{j}^{\prime }.
\]%
Therefore, the conditions are satisfied by the mapping
\[
S\rightarrow V_{j}SV_{j}.
\]

\end{proof}

\begin{myLemma}
\label{lem:IneqGenDCAR} For $T\in M_{n}^{\prime }$ and $E,K_{E}$ as in
theorem \ref{theorem:main},
\[
\Big\Vert\mathcal{D}(T)\Big\Vert_{C_{E}}\leq K_{E}\Big\Vert\left\vert \nabla
_{s}T\right\vert \Big\Vert_{C_{E}}.
\]
\end{myLemma}

\begin{proof}
By lemma \ref{lem:DInvSignsCAR},
\[
\Big\Vert\mathcal{D}(T)\Big\Vert_{C_{E}}=\mathbb{E}\Big\Vert%
\sum_{j=1}^{n}\varepsilon _{j}P_{j}^{\prime }D_{j}^{\prime }(T)\Big\Vert%
_{C_{E}}.
\]

Again we consider two cases for Khintchine inequalities.

\paragraph*{Case 1:}
E is 2-convex and q-concave:

Owing to the commutation relations, if $A\subset \{1,..,n\}$ and $1\leq
j\leq n,$
\[
P_{j}^{\prime }Q_{A}^{\prime }P_{j}^{\prime }{}=(-1)^{\left\vert
A\right\vert }Q_{A}^{\prime }=P_{1}^{\prime }Q_{A}^{\prime }P_{1}^{\prime
}{}.
\]

In particular if $S_{j}\in M_{n}^{\prime },$
\begin{equation}
\Big\Vert (\sum_{j=1}^{n}P_{j}^{\prime }S_{j}S_{j}^{\ast }P_{j}^{\prime
}{})^{\frac{1}{2}}\Big\Vert _{C_{E}} = \Big\Vert (\sum_{j=1}^{n}P_{1}^{%
\prime }S_{j}S_{j}^{\ast }P_{1}^{\prime }{})^{\frac{1}{2}} \Big\Vert %
_{C_{E}} = \Big\Vert (\sum_{j=1}^{n}\left\vert S_{j}^{\ast }\right\vert
^{2})^{\frac{1}{2}} \Big\Vert _{C_{E}}.  \label{19}
\end{equation}

Hence the term $\max $ which appears in the Khintchine upper inequality is
now
\[
\max \Big\{\Big\Vert(\sum_{j=1}^{n}\left\vert D_{j}^{\prime }(T)^{\ast
}\right\vert ^{2})^{\frac{1}{2}}\Big\Vert_{C_{E}},\Big\Vert%
(\sum_{j=1}^{n}\left\vert D_{j}^{\prime }(T)\right\vert ^{2})^{\frac{1}{2}}%
\Big\Vert_{C_{E}}\Big\}\leq \Big\Vert\left\vert \nabla _{s}T\right\vert %
\Big\Vert_{C_{E}}.
\]

\paragraph*{Case 2:}

E is 2-concave:

By the easy part of non-commutative Khintchine inequality, we have:
\[
\mathbb{E}\Big\Vert\sum_{j=1}^{n}\varepsilon _{j}P_{j}^{\prime
}D_{j}^{\prime }(T)\Big\Vert_{C_{E}}\leq \Big\Vert(\sum_{j=1}^{n}\left\vert
D_{j}^{\prime }(T)\right\vert ^{2})^{\frac{1}{2}}\Big\Vert_{C_{E}}\leq %
\Big\Vert\left\vert \nabla _{s}T\right\vert \Big\Vert_{C_{E}}.
\]
\end{proof}

\begin{proofof}{Theorem \ref{theorem:mainCAR}}
The proof is similar to that of Theorem \ref{theorem:main}, using
the previous lemmas.
\end{proofof}

\subsection{Applications}

Here is the analogue of corollary \ref{cor:exp}, with the same proof.

\begin{myCorollary}
\label{cor:expCAR} Let $E$ be as in theorem \ref{theorem:main} and let $T\in
M_{n}^{\prime }$. Then
\[
\;\frac{1}{2}\tau _{n}(\exp \left\vert T-\tau _{n}(T)Id\right\vert )\leq
\tau _{n}(\exp \frac{\pi ^{2}}{16}\left\vert \nabla _{s}T\right\vert
^{2})\leq \exp \frac{\pi ^{2}}{16}\Vert \left\vert \nabla _{s}T\right\vert
\Vert _{C_{\infty }}^{2}).
\]%
and for $0<\alpha <\frac{1}{2}$,
\[
\;\frac{1}{2}\tau _{n}(\exp \left\vert N^{\prime \alpha }\,T\right\vert
)\leq \tau _{n}(\exp \frac{K_{\alpha }^{2}}{2}\,\left\vert \nabla
_{s}(\,T)\right\vert ^{2})\leq \exp \frac{K_{\alpha }^{2}}{2}\,\Vert
\left\vert \nabla _{s}(T)\right\vert \Vert _{C_{\infty }}^{2}).
\]
\end{myCorollary}

\begin{myCorollary}
\label{cor:concCAR} Let $T\in M_{n}^{\prime }$.

\begin{enumerate}
\item Then
\[
\;\tau _{n}(\int_{t}^{\infty }dF(s))\leq 2\exp (-\frac{4t^{2}}{\pi ^{2}\Vert
\,\left\vert \nabla _{s}T\right\vert \,\Vert _{C_{\infty }}^{2}})
\]%
where $dF$ denotes the spectral measure of $\left\vert T-\tau
_{n}(T)Id\right\vert $.

\item Let $0<\alpha <\frac{1}{2}$. Then
\[
\;\tau _{n}(\int_{t}^{\infty }dG_{\alpha }(s))\leq 2\exp (-\frac{t^{2}}{%
2K_{\alpha }^{2}\Vert \,\left\vert \nabla _{s}(T)\right\vert \,\Vert
_{C_{\infty }}^{2}})
\]

where $dG_{\alpha }$ denotes the spectral measure of $\left\vert N^{\prime
\alpha }\,T\right\vert .$
\end{enumerate}
\end{myCorollary}

\begin{proof}
These inequalities are proved as in corollary \ref{cor:conc}, replacing the
classical Tchebychev inequality e.g. in the first case by (for $t\geq 0)$
\[
e^{t}\tau _{n}(\int_{t}^{\infty }dF(s))\leq \tau _{n}(\int_{t}^{\infty
}e^{s}dF(s))\leq \tau _{n}(\exp \left\vert T-\tau _{n}(T)Id\right\vert ).
\]
\end{proof}

\subsection{The reverse inequality}

The next theorem is the analogue of theorem \ref{5.1}; it was proved in \cite%
{LP} \ for $\beta =\frac{1}{2}$ and $E=L^{p}(\Omega _{n}).$ $K_{E}$ and $%
H_{E}$ are defined as in theorem \ref{theorem:main}.

\begin{myTheorem}
\label{6.8}Let $E$ be a symmetric function space on $\Omega _{n},$ let $%
k_{\beta }$ be defined as in theorem \ref{5.1} for $\beta >\frac{1}{2}.$
Then for $T$ $\ $in the CAR algebra $M_{n}^{\prime }$

\begin{enumerate}
\item if $E$ is 2-convex%
\[
\Big\Vert\left\vert \nabla _{s}T\right\vert \Big\Vert_{C_{E}}\leq 2k_{\beta }%
\Big\Vert N^{\prime \beta }(T)\Big\Vert_{C_{E}}
\]

\item if $E$ is 2-concave and $r-$convex%
\[
\inf_{D_{j}^{\prime }(T)=V_{j}+W_{j}} \Big\{ \left\Vert (\sum \left\vert
V_{j}\right\vert ^{2})^{\frac{1}{2}}\right\Vert _{C_{E}}+\left\Vert (\sum
\left\vert W_{j}^{\ast }\right\vert ^{2})^{\frac{1}{2}}\right\Vert _{C_{E}} %
\Big\} \leq k_{\beta }K_{E^{\ast }}^{2}\Big\Vert N^{\prime \beta }(T)%
\Big\Vert_{C_{E}}.
\]

\item Moreover, if $C_{E}$ is UMD, similar inequalities hold for $\beta =%
\frac{1}{2},$ replacing $k_{\beta }$ by $H_{E}.$
\end{enumerate}
\end{myTheorem}

\noindent Before proving this theorem we first adapt the notation and lemmas
of part 5.

Let $\Pi _{j}^{\prime }$ be the orthogonal projection from $C_{2}(\mathcal{M}%
_{2^{n}})$ onto $P_{j}^{\prime }C_{2}(M_{n}^{\prime }),$ $1\leq j\leq n,$
and let $\Pi ^{\prime }=\sum\limits_{j=1}^{n}\Pi _{j}^{\prime }.$

\noindent The analogue of lemma \ref{5.4} 2.(i) is still valid for $\Pi
^{\prime }.$ Lemmas \ref{5.3} and \ref{5.4} 2.(ii) are replaced by the
following:

\begin{myLemma}
\label{6.9}Let $E$ be a symmetric sequence space on $\Omega _{n}.$

\begin{enumerate}
\item Let $E$ be 2-convex.

Then for $S\in \mathcal{M}_{2^{n}}$ and $S_{j}=P_{j}^{\prime }\Pi
_{j}^{\prime }(S)\in M_{n}^{\prime },$ $1\leq j\leq n,$%
\[
\max \Big\{ \left\Vert (\sum \left\vert S_{j}\right\vert ^{2})^{\frac{1}{2}%
}\right\Vert _{C_{E}(M^{\prime}_{n})},\left\Vert (\sum \left\vert
S_{j}^{\ast }\right\vert ^{2})^{\frac{1}{2}}\right\Vert
_{C_{E}(M^{\prime}_{n})} \Big\} \leq \left\Vert S\right\Vert _{C_{E}(%
\mathcal{M}_{2^{n}})}.
\]

\item Let $E$ be 2-concave and $r$-convex ($r>1).$

Then for $T_{j}\in C_{E}(M^{\prime}_{n})$ and decompositions $%
T_{j}=V_{j}+W_{j}$ in $C_{E}(M^{\prime}_{n}),$
\[
\inf_{T_{j}=V_{j}+W_{j}} \Big\{ \left\Vert (\sum \left\vert V_{j}\right\vert
^{2})^{\frac{1}{2}}\right\Vert _{C_{E(M^{\prime}_{n})}}+\left\Vert (\sum
\left\vert W_{j}^{\ast }\right\vert ^{2})^{\frac{1}{2}}\right\Vert
_{C_{E}(M^{\prime}_{n})} \Big\} \leq K_{E^{\ast }}\left\Vert \sum
P_{j}^{\prime }T_{j}\right\Vert _{C_{E}(\mathcal{M}_{2^{n}})}.
\]
\end{enumerate}
\end{myLemma}

\begin{proof}
We repeat the proofs of lemmas \ref{5.3},  \ref{5.4} and take (\ref{19})
into account in order to replace in the formulas $|S_{j}^{\ast
}P_{j}^{\prime }|$, $|W_{j}^{\ast }P_{j}^{\prime }|$ by $|S_{j}^{\ast }|$, $%
|W_{j}^{\ast }|.$
\end{proof}

The analogue of lemma \ref{5.2} is still valid:

\begin{myLemma}
\label{6.10}Let $T \in M^{\prime}_{n}$ such that $\tau (T)=0.$ Then, for $%
\beta > 0,$
\begin{equation}
\Gamma (\beta) D^{\prime}_{j}(T)= P^{\prime}_{j} \Pi^{\prime}_{j}
\int_{0}^{ \frac{\pi }{2} } e^{\theta \mathcal{D}} (N'^{\beta
}T)(-Log\cos \theta )^{\beta -1}d\theta ,\;1\leq j\leq n,
\label{***}
\end{equation}
whence
\begin{equation}
\Gamma (\beta )\sum P^{\prime}_{j}D^{\prime}_{j}(T)=\Pi^{\prime}\int_{0}^{%
\frac{\pi }{2}}e^{\theta \mathcal{D}}(N'^{\beta}T)(-Log\cos \theta
)^{\beta -1}d\theta .  \label{23}
\end{equation}
\end{myLemma}

\begin{proof}
If $j\in A,$ in virtue of the commutation relations,
\[
P_{j}^{\prime }\Pi _{j}^{\prime }e^{\theta \mathcal{D}}(Q_{A}^{\prime
})=\sin \theta \cos ^{\left\vert A\right\vert -1}\theta P_{j}^{\prime
}(\prod\limits_{k\in A,k<j}\,Q_{k}^{\prime })P_{j}^{\prime
}\prod\limits_{k\in A,k>j}Q_{k}^{\prime }
\]
\[
=\sin \theta \cos ^{\left\vert A\right\vert -1}\theta Q_{j}^{\prime
}(\prod\limits_{k\in A,k<j}\,Q_{k}^{\prime })\,Q_{j}^{\prime
}\prod\limits_{k\in A,k>j}Q_{k}^{\prime }=\sin \theta \cos ^{\left\vert
A\right\vert -1}\theta D_{j}^{\prime }(Q_{A}^{\prime }).
\]
\noindent By (\ref{31}) this implies (\ref{***}) hence (\ref{23}).
\end{proof}

\begin{proofof}{Theorem \ref{6.8}}

\paragraph*{Case 1:}

$E$ is 2-convex.

Applying lemma \ref{6.9}, case 1, to $S_{j}=D_{j}^{\prime }(T)$, we get by (%
\ref{***})
\[
\frac{1}{2}\Big\Vert\left\vert \nabla _{s}T\right\vert \Big\Vert_{C_{E}}\leq
\max \Big\{\Big\Vert(\sum_{j=1}^{n}\left\vert D_{j}^{\prime }(T)^{\ast
}\right\vert ^{2})^{\frac{1}{2}}\Big\Vert_{C_{E}},\Big\Vert%
(\sum_{j=1}^{n}\left\vert D_{j}^{\prime }(T)\right\vert ^{2})^{\frac{1}{2}}%
\Big\Vert_{C_{E}}\Big\}\leq k_{\beta }\Big\Vert N^{\prime \beta }(T)\Big\Vert%
_{C_{E}}.
\]

\paragraph*{Case 2:}

$E$ is 2-concave and $r-$convex.

\noindent By lemma \ref{6.9}, case 2, then by the analogue of lemma \ref{5.4}
2.(i) and (\ref{23}), we get%
\[
\inf_{D_{j}^{\prime }(T)=V_{j}+W_{j}}\Big\{\left\Vert (\sum \left\vert
V_{j}\right\vert ^{2})^{\frac{1}{2}}\right\Vert _{C_{E}}+\left\Vert (\sum
\left\vert W_{j}^{\ast }\right\vert ^{2})^{\frac{1}{2}}\right\Vert _{C_{E}}%
\Big\}\leq K_{E^{\ast }}\left\Vert \sum P_{j}^{\prime }D_{j}^{\prime
}(T)\right\Vert _{C_{E}}
\]%
\[
\leq k_{\beta }K_{E^{\ast }}^{2}\left\Vert N^{\prime \beta }(T)\right\Vert
_{C_{E}}.
\]

\paragraph*{The UMD case:}

This case is a modification of cases 1 and 2 exactly as in the proof of
theorem \ref{5.1}.
\end{proofof}

\appendix

\section{\textbf{A refinement of theorem \protect\ref{theorem:main}}}

All the results in this paper can be achieved activating the
derivative and the number operator only on part of the coordinates.
We give an example for theorem \textbf{\ref{theorem:main}}case
1\textbf{. }The proof is the same.

\noindent Let $\mathcal{J}\subset \{1,..,n\}$ be an arbitrary set of
coordinates, and let $\bar{\mathcal{J}}=\{1,..,n\}\setminus \mathcal{J}$. We
associate to $f:\Omega _{n}\rightarrow \mathbb{C}$ the following functions

\[
\mathcal{V}_{\mathcal{J}}(f)=\sum_{A\cap \mathcal{J}\neq \emptyset }\hat{f}%
(A)\omega _{A},\;\mathcal{P}_{\mathcal{J}}(f)=\sum_{A\subseteq \mathcal{J}}%
\hat{f}(A)\omega _{A}.
\]

\noindent so that $f=\mathcal{V}_{\mathcal{J}}(f)+\mathcal{P}_{\bar{\mathcal{%
J}}}(f)$. Let $|\nabla _{\mathcal{J}}(f)|=\sum_{j\in \mathcal{J}}(|\partial
_{j}(f)|^{2})^{\frac{1}{2}}$. Let

\[
\mathcal{R}_{\theta ,\mathcal{J}}=R_{\theta }^{\delta _{\mathcal{J}
}(1)}\otimes ...\otimes R_{\theta }^{\delta _{\mathcal{J}}(n)}\in \mathcal{M}%
_{2^{n}},\ \ \theta \in \mathbb{R},
\]

\noindent where $R_{\theta }^{\delta _{\mathcal{J}}(j)}$ occurs in the j$%
^{th}$ factor, and
\[
R_{\theta }^{\delta _{\mathcal{J}}(j)}=\left\{
\begin{array}{ll}
Id & \mbox{if $j \not\in \calJ$} \\
R_{\theta } & \mbox{if $j \in \calJ$}%
\end{array}
\right. .
\]

\noindent The action $T\rightarrow \mathcal{R}_{\theta ,\mathcal{J}}^{\ast
}T \mathcal{R}_{\theta ,\mathcal{J}}$ is again an inner *-automorphism of $%
\mathcal{M}_{2^{n}}$ which preserves the trace $\tau _{n}$ and as in lemma %
\ref{lem:factFormula} (where $\mathcal{J=}\{1,..,n\})$ we have

\[
\mathcal{R}_{\theta ,\mathcal{J}}^{\ast }Q_{A}\mathcal{R}_{\theta ,\mathcal{J%
}}=\prod_{j\in A\cap \bar{\mathcal{J}}}Q_{j}\prod\limits_{j\in A\cap
\mathcal{J}}(\cos \theta \,Q_{j}+\sin \theta \,P_{j})=e^{\theta \mathcal{D}_{%
\mathcal{J}}}(T).
\]

\noindent For $T\in M_{n}$, let $N_{\mathcal{J}}=\sum_{j\in \mathcal{J}%
}D_{j}^{\ast }(T)D_{j}(T)$. The analogue of (\ref{equ:factorFormula}) is
\[
\mathcal{E}_{M_{n}}e^{\theta \mathcal{D}_{\mathcal{J}}}J_{n}(T)=\cos ^{N_{%
\mathcal{J}}}\theta (T).
\]

\begin{myTheorem}
\label{theo:mainPartialSetOfCoordinates} Let $E$ be a symmetric function
space on $\Omega _{n}$. Then for every function $f:\Omega _{n}\rightarrow
\mathbb{C}$ and $\mathcal{J}\subset \{1,..,n\}$ we have
\[
\Big\Vert f-\mathcal{P}_{\bar{\mathcal{J}}}(f)\Big\Vert_{E}=\Big\Vert%
\mathcal{V}_{\mathcal{J}}(f)\Big\Vert_{E}\leq \frac{\pi }{4}K_{E}\Big\Vert%
|\nabla _{\mathcal{J}}(f)|\Big\Vert_{E}.
\]
\end{myTheorem}

\begin{myCorollary}
Let $\mathcal{J}\subset \{1,..,n\}$. Then for every $f:\Omega
_{n}\rightarrow \mathbb{C}$,

\begin{enumerate}
\item
\[
\Vert f\Vert _{p}\leq C\sqrt{p}\,\Vert \,|\nabla _{\mathcal{J}}(f)|\,\Vert
_{p}+\Vert \mathcal{P}_{\bar{\mathcal{J}}}(f)\Vert _{p}\;\;1\leq p<\infty
\]

where $C$ is a universal constant. In particular if $\bar{\mathcal{J}}=\{j\}$
\[
\Vert f\Vert _{p}\leq C\sqrt{p}\,\Vert \,|\nabla _{\{1,..,n\}\setminus
j}(f)|\,\Vert _{p}+|\hat{f}(\{j\})|.
\]

\item
\[
\ \left\Vert f-\mathcal{P}_{\bar{\mathcal{J}}}(f)\right\Vert _{\infty }\leq C%
\sqrt{n}\,\Vert \,|\nabla _{\mathcal{J}}(f)|\,\Vert _{\infty }
\]
\end{enumerate}
\end{myCorollary}

\begin{proof}
Indeed, for $f:\Omega _{n}\rightarrow \mathbb{C}$, $\ \frac{1}{2}\Vert
f\Vert _{L^{\infty }(\Omega _{n})}\leq \Vert f\Vert _{L^{n}(\Omega
_{n})}\leq \Vert f\Vert _{L^{\infty }(\Omega _{n})}$. Hence the last
inequality \ follows from the first which is a particular case of theorem %
\ref{theo:mainPartialSetOfCoordinates}.
\end{proof}

\end{document}